\documentclass[leqno]{siamltex}

\usepackage[letterpaper=true, colorlinks=true, linkcolor=red, filecolor=green, citecolor=red, pdfpagemode=None]{hyperref}

\usepackage{amsmath}
\usepackage{amssymb}
\usepackage{graphicx}
\usepackage{epstopdf}
\usepackage{color}
\usepackage{algorithm}

\newtheorem{example}{Example}[section]
\newcommand{\atopfrac}[2]{\genfrac{}{}{0pt}{}{#1}{#2}}

\DeclareMathOperator{\tridiag}{tridiag}

\def\qed{~\vbox{\hrule\hbox{\vrule height1.3ex\hskip0.8ex\vrule}\hrule}}

\newenvironment{dedication}
        {\vspace{1em}\begin{quotation}\begin{center}\begin{em}}
        {\par\end{em}\end{center}\end{quotation}\vspace{1.5em}}

\title{Block diagonal dominance of matrices revisited: bounds for the norms of inverses and\\ eigenvalue inclusion sets\thanks{This work dated May 16, 2018.}}

\author{Carlos Echeverr{\'i}a\footnotemark[4] , J{\"o}rg Liesen\footnotemark[4] ,\and Reinhard Nabben\footnotemark[4]}

\begin{document}

\maketitle

\renewcommand{\thefootnote}{\fnsymbol{footnote}}

\footnotetext[4]{{Institute of Mathematics, Technical University of Berlin, Strasse des 17. Juni 136, D-10623 Berlin, Germany. Carlos Echeverria's work was partially supported by the Einstein Center for Mathematics, Berlin, the Deutsche Akademische Austausch Dienst (DAAD), Germany,  and the Consejo Nacional de Ciencia y Tecnolog\'ia (CONACyT), M{\'e}xico.\\ Contact e-mails: {\tt \{echeverria, liesen, nabben\}@math.tu-berlin.de}}. }

\begin{dedication}
{\scriptsize Dedicated to Richard S. Varga}
\end{dedication}

\begin{abstract}
We generalize the bounds on the inverses of diagonally dominant matrices obtained in~\cite{Nabben99} from scalar to block tridiagonal matrices.
Our derivations are based on a generalization of the classical condition of block diagonal dominance of matrices given by Feingold and Varga in~\cite{FeiVar62}. Based on this generalization, which was recently presented in \cite{BenEvaHamLupSla17}, we also derive a variant of the Gershgorin Circle Theorem for general block matrices which can provide tighter spectral inclusion regions than those obtained by Feingold and Varga.
\end{abstract}



\pagestyle{myheadings}
\thispagestyle{plain}
\markboth{C. ECHEVERR{\'I}A, J. LIESEN and R. NABBEN}{BLOCK DIAGONAL DOMINACE AND ITS APPLICATIONS}

\section{Introduction}\label{s:intro}

Matrices that are characterized by off-diagonal decay, or more generally ``localization'' of their {entries}, appear in applications throughout the mathematical and computational sciences. The presence of such localization can lead to  computational savings, since it allows to (closely) approximate a given matrix by using its significant entries only, and discarding the negligible ones according to a pre-established criterion. In this context it is then of great practical interest to know a priori how many and which of these entries can be discarded as insignificant. Many authors have therefore studied decay rates for different matrix classes and functions of matrices; see, e.g.,~\cite{BenBoi14,BenGol99,BenRaz07,BenSim15,CanSimVer14,DemMosSmi84,KriStrWer15,PelPol01}. For an excellent survey of the current state-of-the-art we refer to~\cite{Ben16}.

An important example in this context is given by the (nonsymmetric) diagonally dominant matrices, and in particular the diagonally dominant tridiagonal matrices, which were studied, e.g., in~\cite{Nabben99-2,Nabben99}. As shown in these works, the entries of the inverse decay with an exponential rate along a row or column, depending on whether the given matrix is row or column diagonally dominant; see~\cite[Section~3.2]{Ben16} for a more general treatment of decay bounds for the inverse and further references. Our main goal in this paper is to generalize results of~\cite{Nabben99} from scalar to block tridiagonal matrices. In order to do so, we use a generalization of the classical definition of block diagonal dominance
of Feingold and Varga~\cite{FeiVar62} to derive bounds
and decay rates for the block norms of the inverse of block tridiagonal matrices. We also show how to improve these bounds iteratively (Section~\ref{s:blkmat}).
Moreover, we obtain a new variant of the Gershgorin Circle Theorem for general block matrices (Section~\ref{s:eigenregion}). Throughout this paper we assume that $\|\cdot\|$ is a given submultiplicative matrix norm.

\section{Bounds on the inverses of block tridiagonal matrices}\label{s:blkmat}
We start with a definition of block diagonally dominant matrices, which was recently presented in \cite{BenEvaHamLupSla17} in an application of block diagonal preconditioning.

\begin{definition}\label{def:bdd}
Consider a matrix of the form
\begin{equation}\label{eq:blockmat}
A=[A_{ij}]\quad
\mbox{with blocks $A_{ij}\in\mathbb{C}^{m\times m}$ for $i,j=1,\dots,n$.}
\end{equation}
The matrix $A$ is called \emph{row block diagonally dominant} (with respect to the matrix norm $\|\cdot\|$)
when the diagonal blocks $A_{ii}$ are nonsingular, and

\begin{equation}\label{eq:blocdiagdom1}
\sum_{\atopfrac{j=1}{j\neq i}}^{n}\|A_{ii}^{-1}A_{ij}\|\leq 1,\quad\text{for $i=1,\dots,n$.}
\end{equation}
If a strict inequality holds in \eqref{eq:blocdiagdom1}
then $A$ is called \emph{row block strictly diagonally dominant} (with respect to the matrix norm $\|\cdot\|$).
\end{definition}

Obviously, an analogous definition of \emph{column} block diagonal dominance is possible. Most of the results in this paper can be easily rewritten for that case. Also note that the authors of~\cite{BenEvaHamLupSla17} call a matrix block diagonally dominant when all its diagonal
blocks are nonsingular, and \eqref{eq:blocdiagdom1} or the anologous conditions with $A_{ij}A_{ii}^{-1}$ replacing $A_{ii}^{-1}A_{ij}$
hold (in the $1$-norm).

The above definition of (row) block diagonal
dominance generalizes the one of Feingold and Varga given in \cite[Definition~1]{FeiVar62}, who considered a
matrix as in \eqref{eq:blockmat} block diagonally dominant when the diagonal blocks $A_{ii}$ are nonsingular,
and
\begin{equation}\label{eq:blocdiagdomFV}
\sum_{\atopfrac{j=1}{j\neq i}}^{n}\|A_{ij}\|\leq(\|A_{ii}^{-1}\|)^{-1},\quad\text{for $i=1,\dots,n$}.
\end{equation}
It is clear that if a matrix satisfies these conditions, then it also satisfies the conditions given in Definition~\ref{def:bdd}. According to Varga~\cite[p.~156]{Var2011}, the definition of block diagonal dominance given in~\cite{FeiVar62} is one of the earliest, and it was roughly simultaneously and independently considered also by Ostrowski~\cite{Ost61} and Fiedler and Pt\'ak~\cite{FiePta62}. Varga calls this a ``Zeitgeist'' phenomenon.

{In} the special case $m=1$, i.e., all the blocks $A_{ij}$ are of size $1\times 1$ and $\|{A_{ij}}\|=|A_{ij}|$, the inequalities \eqref{eq:blocdiagdom1} 
and \eqref{eq:blocdiagdomFV} are equivalent, and they can all be written as
\[
\sum_{\atopfrac{j=1}{j\neq i}}^{n}|A_{ij}|\leq |A_{ii}|,\quad\text{for $i=1,\dots,n$,}
\]
which is the usual definition of row diagonal dominance.

In the rest of this section we will restrict our attention to block tridiagonal matrices of the form
\begin{equation}\label{eq:blocktridiag}
A=\left[ \begin{array}{crrrr}
A_1	& 	B_1		& 					&					&\\
C_1	&	A_2		& B_2			&					&\\
		& \ddots & \ddots 		& \ddots 		& \\
		&				& C_{n-2} 	& A_{n-1}	& B_{n-1}\\
		&				&					& C_{n-1}		& A_{n}
\end{array} \right],\quad \mbox{where $A_i,B_i,C_i\in {\mathbb C}^{m\times m}$.}
\end{equation}
First Capovani for the scalar case in \cite{Cap70, Cap71} and later Ikebe for the block case in~\cite{Ike79} (see also~\cite{Nabben99}), have shown that the inverse of a nonsingular block tridiagonal
matrix can be described by four sets of matrices. The main result can be stated as follows.
\begin{theorem}\label{teo:Ikebe}
Let $A$ be as in \eqref{eq:blocktridiag}, and suppose that $A^{-1}$ as well as $B_i^{-1}$ and $C_i^{-1}$ for
$i=1,\dots,n-1$ exist. If we write $A^{-1}=[Z_{ij}]$ with $Z_{ij}\in {\mathbb C}^{m\times m}$, then there
exist matrices $U_i,V_i,X_i,Y_i\in {\mathbb C}^{m\times m}$ with $U_iV_i=X_iY_i$
for $i=1,\dots,n$, and
\begin{equation}\label{eq:inverse}
Z_{ij}=\begin{cases}
U_{i}V_{j}
&\mbox{if  } i\leq j,
\\
Y_{i}X_{j}
&\mbox{if  } i\geq j.
\end{cases}
\end{equation}
Moreover, the matrices $U_i,V_i,X_i,Y_i$, $i=1,\dots,n$, are recursively given by
\begin{align}
U_1 & = I, \quad U_2=-B_1^{-1}A_1U_1, \label{def:U2}\\
U_i & =-B_{i-1}^{-1}(C_{i-2}U_{i-2}+A_{i-1}U_{i-1}), \quad \mbox{for $i=3,\dots,n$,} \label{def:Ui} \\
V_n & =(A_nU_n+C_{n-1}U_{n-1})^{-1}, \quad
V_{n-1} = -V_nA_nB_{n-1}^{-1}, \label{def:Vn-1}\\
V_{i} &= -(V_{i+1}A_{i+1}+V_{i+2}C_{i+1})B_{i}^{-1},\quad \mbox{for $i=n-2,\dots,1$.} \label{def:Vi}\\
X_1 &= I,\quad X_2=-X_1A_1C_1^{-1}, \\
X_i &= -(X_{i-2}B_{i-2}+X_{i-1}A_{i-1})C_{i-1}^{-1},\quad \mbox{for $i=3,\dots,n$,} \\
Y_n &= (X_nA_n+X_{n-1}B_{n-1})^{-1},\quad Y_{n-1}=-C_{n-1}^{-1}A_nY_n,  \label{def:Yn-1}\\
Y_i &= -C_i^{-1}(A_{i+1}Y_{i+1}+B_{i+1}Y_{i+2}),\quad \mbox{for $i=n-2,\dots,1$.} \label{def:Yi}
\end{align}
\end{theorem}

The next result is a generalization of~\cite[Theorem~3.2]{Nabben99-2}.

\begin{lemma}\label{lem:sequences}
Let $A$ be a matrix as in Theorem~\ref{teo:Ikebe}. Suppose in addition that $A$ is row block diagonally dominant, and that
\begin{equation}\label{eq:blktridom_b}
\|A_1^{-1}B_1\|<1 \quad\mbox{and}\quad \|A_n^{-1}C_{n-1}\|<1.
\end{equation}
Then the sequence $\{\|U_i\|\}_{i=1}^n$ is strictly increasing, and the sequence $\{\|Y_i\|\}_{i=1}^n$ is
strictly decreasing.
\end{lemma}

\begin{proof}
First we consider the
sequence $\{\|U_i\|\}_{i=1}^n$. The definition of $U_2$ in \eqref{def:U2} implies that
$U_1=-A_1^{-1}B_1U_2$. Taking norms and using the first inequality in \eqref{eq:blktridom_b} yields
\[
\|U_1\|\leq\|A_1^{-1}B_1\|\|U_2\|<\|U_2\|.
\]
Now suppose that $\|U_1\|<\|U_2\|<\cdots <\|U_{i-1}\|$ holds for some $i\geq 3$. The equation for $U_i$
in \eqref{def:Ui} can be written as
\[
 -A^{-1}_{i-1}B_{i-1}U_i=U_{i-1}+A_{i-1}^{-1}C_{i-2}U_{i-2}.
\]
 Rearranging terms and taking norms we obtain
 \begin{align*}
 \|U_{i-1}\| &\leq\|A_{i-1}^{-1}C_{i-2}\|\|U_{i-2}\|+\|A_{i-1}^{-1}B_{i-1}\|\|U_{i}\|\\
&<\|A_{i-1}^{-1}C_{i-2}\|\|U_{i-1}\|+\|A_{i-1}^{-1}B_{i-1}\|\|U_{i}\|,
\end{align*}
where we have used the induction hypothesis, i.e., $\|U_{i-2}\|<\|U_{i-1}\|$, in order to obtain the
strict inequality. Since $A$ is row block diagonally dominant we have
\[\|A_{i-1}^{-1}B_{i-1}\|+ \|A_{i-1}^{-1}C_{i-2}\|\leq 1.\]
Combining this with the previous inequality gives
\[
\frac{\|U_{i-1}\|}{\|U_{i}\|}<\frac{\|A_{i-1}^{-1}B_{i-1}\|}{1-\|A_{i-1}^{-1}C_{i-2}\|}\leq 1,
\]
so that indeed $\|U_{i-1}\|<\|U_i\|$.

Next we consider the sequence $\{\|Y_i\|\}_{i=1}^n$. The definition of $Y_{n-1}$ in \eqref{def:Yn-1} implies that $-Y_{n}=A^{-1}_{n}C_{n-1}Y_{n-1}$. Taking norms and using the second inequality in \eqref{eq:blktridom_b} yields
\[
\|Y_{n}\|\leq\|A_n^{-1}C_{n-1}\|\|Y_{n-1}\| < \|Y_{n-1}\|.
\]
Now suppose that $\|Y_n\|<\|Y_{n-1}\|<\cdots<\|Y_{i+1}\|$ holds for some $i\leq n-2$. The equation for
$Y_i$ in \eqref{def:Yi} can be written as
\[
-A_{i+1}^{-1}C_iY_{i}=Y_{i+1}+A_{i+1}^{-1}B_{i+1}Y_{i+2}.
\]
Rearranging terms and taking norms we obtain
\begin{align*}
\|Y_{i+1}\| &\leq \|A_{i+1}^{-1}C_{i}\|\|Y_{i}\|+\|A_{i+1}^{-1}B_{i+1}\|\|Y_{i+2}\| \\
&< \|A_{i+1}^{-1}C_{i}\|\|Y_{i}\|+\|A_{i+1}^{-1}B_{i+1}\|\|Y_{i+1}\|,
\end{align*}
where we have used the induction hypothesis, i.e., $\|Y_{i+2}\|<\|Y_{i+1}\|$, in order to obtain the
strict inequality. Since $A$ is row block diagonally dominant we have
\[
\|A_{i+1}^{-1}C_{i}\|+\|A_{i+1}^{-1}B_{i+1}\|\leq 1.
\]
Combining this with the previous inequality gives
\begin{equation*}
\frac{\|Y_{i+1}\|}{\|Y_{i}\|}<\frac{\|A_{i+1}^{-1}C_{i}\|}{1-\|A_{i+1}^{-1}B_{i+1}\|}\leq 1
\end{equation*}
so that indeed $\|Y_{i+1}\| <\|Y_i\|$.
\end{proof}

For the rest of this section we will assume that that $A$ is a matrix as in Lemma~\ref{lem:sequences}.
Then the inverse is given by $A^{-1}=[Z_{ij}]$
with $Z_{ij}=Y_iX_j$ for $i\geq j$; see Theorem~\ref{teo:Ikebe}. Thus, for each \emph{fixed}
$j=1,\dots, n$, the strict decrease of the sequence $\{\|Y_i\|\}_{i=1}^n$ suggests that
the sequence $\{\|Z_{ij}\|\}_{i=j}^n$ decreases as well, i.e., that the norms of the
blocks of $A^{-1}$ decay columnwise away from the diagonal. We will now study this decay
in detail.

\medskip
We set $C_0=B_n=0$, and define
\begin{align*}
\tau_i   &= \frac{\|A_i^{-1}B_i\|}{1-\|A_i^{-1}C_{i-1}\|},\quad \text{for $i=1,\ldots,n$}, \\
\omega_i &= \frac{\|A_i^{-1}C_{i-1}\|}{1-\|A_i^{-1}B_{i}\|},\quad \text{for $i=1,\ldots,n$}.
\end{align*}

\noindent The row block diagonal dominance of $A$ then implies that $0\leq \tau_i\leq 1$ and $0~\leq~\omega_i~\leq~1$. Also note that, by assumption, $\tau_1=\|A_1^{-1}B_1\|<1$,\linebreak $\omega_n=~\|A_n^{-1}C_{n-1}\|~<~1$, and $\tau_n=\omega_1=0$.

In order to obtain bounds on the norms of the block entries $A^{-1}$, we will first derive alternative recurrence formulas for the matrices $U_i$ and $Y_i$ from Lemma~\ref{lem:sequences}. To this end, we introduce some intermediate quantities and give bounds on their norms in the following result.

\newpage
\begin{lemma}\label{l:nonsing}
The following assertions hold:

\begin{itemize}
\item[(a)]
The matrices $L_1=T_1=A^{-1}_1B_1$, $T_2=I-A_2^{-1}C_1T_1$, and
\begin{align*}
L_i &= T_i^{-1}A_i^{-1}B_{i}, \qquad\qquad \text{for $i=2,\dots,n-1$}, \\
T_i &= I-A_i^{-1}C_{i-1}L_{i-1},\quad  \text{for $i=3,\dots,n-1$,}
\end{align*}
are all nonsingular, and $\|L_i\|\leq\tau_i$, for $i=1,\ldots,n-1$.

\item[(b)]
The matrices $M_n=W_n=A^{-1}_{n}C_{n-1}$, $W_{n-1}=I-A^{-1}_{n-1}B_{n-1}W_{n}$, and
\begin{align*}
M_i&= W_{i}^{-1}A_i^{-1}C_{i-1},\quad\quad \text{for $i=n-1,\dots,2$},\\ 
W_i&= I-A_{i}^{-1}B_{i}M_{i+1}, \quad\text{for $i=n-2,\dots,2$}, 
\end{align*}
are all nonsingular, and $\|M_i\|\leq\omega_i$, for $i=2,\dots,n$.

\end{itemize}
\end{lemma}

\begin{proof}
We only prove \emph{(a)}; the proof of \emph{(b)} is analogous.
The matrices $L_1=T_1=A_1^{-1}B_1$ are nonsingular since both $A_1$ and $B_1$ are. Moreover, \eqref{eq:blktridom_b} gives
$\|L_1\|=\|T_1\|=\|A_1^{-1}B_1\|=\tau_1<1$.
Now suppose that $\|L_{i-1}\| \leq \tau_{i-1}\leq 1$ holds for some $i\geq 2$. Then
\[\|A_i^{-1}C_{i-1}L_{i-1}\|\leq\|A_i^{-1}C_{i-1}\|\|L_{i-1}\|<1,\]
where we have also used that $\|A_i^{-1}C_{i-1}\|\leq 1-\|A_i^{-1}B_i\|<1$. Thus,\linebreak
$T_i=I-A_i^{-1}C_{i-1}L_{i-1}$ is nonsingular, and therefore $L_i=T_i^{-1}A_i^{-1}B_i$
is nonsingular. Using the Neumann series gives
\begin{align*}
\|T_i^{-1}\|&=\|(I-A_i^{-1}C_{i-1}L_{i-1})^{-1}\|=\left\|\sum_{k=0}^{\infty}(A_i^{-1}C_{i-1}L_{i-1})^{k}\right\|
\leq\sum_{k=0}^{\infty}\|A_i^{-1}C_{i-1}L_{i-1}\|^{k}\\
&=\frac{1}{1-\|A_i^{-1}C_{i-1}L_{i-1}\|}\leq\frac{1}{1-\|A_i^{-1}C_{i-1}\|\|L_{i-1}\|}
\leq\frac{1}{1-\|A_i^{-1}C_{i-1}\|},
\end{align*}
and $\|L_i\|=\|T_i^{-1}A_i^{-1}B_i\|\leq\frac{\|A_i^{-1}B_i\|}{1-\|A_i^{-1}C_{i-1}\|}=\tau_i \leq 1$,
which finishes the proof.
\end{proof}

Using Lemma~\ref{l:nonsing} we can now derive alternative recurrences for the matrices $U_i$ and $Y_i$
from Lemma~\ref{lem:sequences}.

\begin{lemma}\label{l:recurrenceUY}
If $A$ is a matrix as in Lemma~\ref{lem:sequences}, then the corresponding matrices $U_i$ and $Y_i$
are given by
\begin{align}
U_i &= -L_iU_{i+1},  \quad\text{for $i=1,\ldots n-1$},\label{eq:Us1}\\
Y_i &= -M_{i}Y_{i-1}, \quad \text{for $i=n,\ldots 2$},\label{eq:Us2}
\end{align}
where the matrices $L_i$ and $M_i$ are defined as in Lemma~\ref{l:nonsing}.
\end{lemma}

\begin{proof}
We only prove that \eqref{eq:Us1} holds; the proof of \eqref{eq:Us2} is analogous.
From \eqref{def:U2} and the definition of $T_1$ in Lemma~\ref{lem:sequences} we obtain
\begin{equation*}
U_1=-A_1^{-1}B_1U_2=-L_1U_2.
\end{equation*}
We next write \eqref{def:Ui} for $i=3$ as
\[
-A_2^{-1}B_2U_3=A_2^{-1}C_1U_1+U_2=-A_2^{-1}C_1T_1U_2+U_2=(I-A_2^{-1}C_1T_1)U_2=T_2U_2,
\]
and hence
\[
U_2=-T_2^{-1}A_2^{-1}B_2U_3=-L_2U_3.
\]
Now suppose that $U_{i-1}=-L_{i-1}U_i$ holds for some $3\leq i\leq n-1$. Then from~\eqref{def:Ui} we obtain
\[
-A_{i}^{-1}B_iU_{i+1}=A_i^{-1}C_{i-1}U_{i-1}+U_{i}=(I-A_i^{-1}C_{i-1}L_{i-1})U_{i}=T_iU_i,
\]
and hence
\[
U_i=-T_i^{-1}A_i^{-1}B_iU_{i+1}=-L_iU_{i+1},
\]
which completes the proof.
\end{proof}

We are now ready to state and prove our bounds on the norms of the blocks of~$A^{-1}$,
which generalize~\cite[Theorems~3.1 and~3.2]{Nabben99} from the scalar to the block case.

\begin{theorem}\label{teo:blockbounds}
If $A$ is a matrix as in Lemma~\ref{lem:sequences}, then $A^{-1}=[Z_{ij}]$ with
\begin{align}
\|Z_{ij}\| & \leq\|Z_{jj}\| \prod_{k=i}^{j-1}\tau_k,\qquad\text{for all $i<j$}, \label{offdiagbound1}\\
\|Z_{ij}\| & \leq\|Z_{jj}\| \prod_{k=j+1}^{i}\omega_k,\quad\text{for all $i>j$}.\label{offdiagbound2}
\end{align}
Moreover, for $i=1,\dots,n$,
\begin{equation}\label{diagbounds}
\frac{\|I\|}{\|A_i\|+\tau_{i-1}\|C_{i-1}\|+\omega_{i+1}\|B_{i}\|}
\leq
\|Z_{ii}\|
\leq
\frac{\|I\|}{\|A_i^{-1}\|^{-1}-\tau_{i-1}\|C_{i-1}\|-\omega_{i+1}\|B_i\|}, 
\end{equation}
provided that the denominator of the upper bound is larger than zero,
and where we set $C_0=B_n=0$, and $\tau_0=\omega_{n+1}=0$.
\end{theorem}

\emph{Proof.}
From Lemma~\ref{l:recurrenceUY} we know that $U_i=-L_iU_{i+1}$ holds for
$i=1,\dots,n~-~1$. Thus, for all $i<j$,
\[Z_{ij}=U_iV_j=-L_iU_{i+1}V_j
=(-1)^{j-i}\left(\prod_{k=i}^{j-1}L_k\right)U_jV_j
=(-1)^{j-i}\left(\prod_{k=i}^{j-1}L_k\right)Z_{jj}.\]
Taking norms and using Lemma~\ref{l:nonsing} yields
\[\|Z_{ij}\|\leq
\|Z_{jj}\|\prod_{k=i}^{j-1}\|L_k\|\leq\|Z_{jj}\|\prod_{k=i}^{j-1}\tau_k.\]
The expression for $i>j$ follows analogously using the two lemmas.

Since $AA^{-1}=AZ=I$ we have
\[
C_{i-1}Z_{i-1,i}+A_iZ_{ii}+B_iZ_{i+1,i}=I,\quad
\mbox{for $i=1,\dots,n$,}
\]
where we set $C_0=Z_{0,1}=B_n=Z_{n+1,n}=0$. Using \eqref{eq:inverse} and Lemma~\ref{l:recurrenceUY},
\begin{align*}
Z_{i-1,i} &=U_{i-1}V_i=-L_{i-1}U_iV_i=-L_{i-1}Z_{ii},\\
Z_{i+1,i} &=Y_{i+1}X_{i}=-M_{i+1}Y_iX_i=-M_{i+1}Z_{ii},
\end{align*}
where we set $U_0=L_0=Y_{n+1}=M_{n+1}=0$. Combining this with the previous equation yields
\begin{equation}\label{eq:inverse1}
-C_{i-1}L_{i-1}Z_{ii}-B_iM_{i+1}Z_{ii}+A_iZ_{ii}=I,
\quad \mbox{for $i=1,\dots,n$,}
\end{equation}
Taking norms and using again Lemma~\ref{l:nonsing} now gives
\begin{align*}
\|I\|&=\|-C_{i-1}L_{i-1}Z_{ii}-B_iM_{i+1}Z_{ii}+A_iZ_{ii}\|\\
&\leq (\|C_{i-1}\|\|L_{i-1}\|+\|A_i\|+\|B_i\|\|M_{i+1}\|)\|Z_{ii}\|\\
&\leq (\tau_{i-1}\|C_{i-1}\|+\|A_i\|+\omega_{i+1}\|B_i\|)\|Z_{ii}\|, \quad\mbox{for $i=1,\dots,n$},
\end{align*}
where we set $\tau_0=\omega_{n+1}=0$, and which shows
the lower bound in \eqref{diagbounds}. In order to show the upper bound we write  \eqref{eq:inverse1} as
\begin{equation*}
I-A_iZ_{ii}= -(C_{i-1}L_{i-1}Z_{ii}+B_iM_{i+1}Z_{ii}),\quad \mbox{for $i=1,\dots,n$.}
\end{equation*}
This yields
\begin{align*}
\|A_iZ_{ii}\|-\|I\| &\leq \|I-A_iZ_{ii}\|
= \|C_{i-1}L_{i-1}Z_{ii}+B_iM_{i+1}Z_{ii}\| \\
 & \leq  (\tau_{i-1}\|C_{i-1}\|+\omega_{i+1}\|B_i\|)\|Z_{ii}\|.
\end{align*}
From
$\|Z_{ii}\|=\|A_i^{-1}A_iZ_{ii}\|\leq\|A_i^{-1}\|\|A_iZ_{ii}\|$
we get $\|A_iZ_{ii}\|\geq \|Z_{ii}\|/\|A_{i}^{-1}\|$,
and combining this with the previous inequality yields
\[\left(\tau_{i-1}\|C_{i-1}\|+\omega_{i+1}\|B_i\|\right)\|Z_{ii}\|\geq\frac{1}{\|A_{i}^{-1}\|}\|Z_{ii}\|-\|I\|.\]
When $\|A_i^{-1}\|^{-1}-\tau_{i-1}\|C_{i-1}\|-\omega_{i+1}\|B_i\|>0$ holds, we get the upper bound
in \eqref{diagbounds}.\qed \\

Note that the positivity assumption on the denominator of the upper bound in \eqref{diagbounds}
is indeed necessary. A simple example for which the denominator
is equal to zero is given by the matrix
$A={\rm tridiag}(-1,2,-1) \in {\mathbb R}^{n\times n}$ with $1\times 1$ blocks,
which satisfies all assumptions of Lemma~\ref{lem:sequences}.

Both the off-diagonal bounds \eqref{offdiagbound1}--\eqref{offdiagbound2} and the diagonal bounds \eqref{diagbounds} depend on the values $\tau_i$ and $\omega_i$, which bound $\|L_{i}\|$ and $\|M_{i}\|$, respectively. We will now show that by modifying the proof of Lemma~\ref{l:nonsing} the
bounds can be improved in an iterative fashion. This is analogous to the iterative improvement
for the case when the blocks of $A$ are scalars, which was considered in~\cite{Nabben99}.

We have shown in the inductive proof of Lemma \ref{l:nonsing} that
\begin{equation*}
\|T_i^{-1}\|\leq\frac{1}{1-\|A_i^{-1}C_{i-1}L_{i-1}\|}\leq\frac{1}{1-\|A_i^{-1}C_{i-1}\|\|L_{i-1}\|}\leq\frac{1}{1-\|A_i^{-1}C_{i-1}\|}. 
\end{equation*}
This bound can be improved by making use of Lemma~\ref{l:nonsing} itself, i.e.,
\begin{equation*}
\|T_i^{-1}\|\leq\frac{1}{1-\|A_i^{-1}C_{i-1}L_{i-1}\|}\leq\frac{1}{1-\|A_i^{-1}C_{i-1}\|\|L_{i-1}\|}\leq\frac{1}{1-\|A_i^{-1}C_{i-1}\|\tau_{i-1}},
\end{equation*}
and this yields
\begin{equation*}
\|L_i\|=\|T_i^{-1}A_i^{-1}B_i\|\leq\frac{\|A_i^{-1}B_i\|}{1-\|A_i^{-1}C_{i-1}\|\tau_{i-1}}.
\end{equation*}
If we denote the expression on the right hand side by $\tau_{i,2}$, then we obtain a modified
version of Lemma~\ref{l:nonsing}, where $\|L_i\|\leq\tau_{i,2}\leq\tau_2\leq 1$. Iteratively
we now define, for all $i=1,\ldots,n$ and $t=1,\ldots,n-1$,
\begin{eqnarray*}
\tau_{i,t}=\begin{cases}
\tau_i&\mbox{if } t=1,\\
\tau_{i,t-1} &\mbox{if } t>i, \\
\frac{\|A_i^{-1}B_i\|}{1-\|A_i^{-1}C_{i-1}\|\tau_{i-1,t-1}} & \mbox{else}.
\end{cases}
\end{eqnarray*}
Analogously we can proceed for the values $\|M_i\|$, and here we define,
for all $i=1,\ldots,n$ and $t=1,\ldots,n-1$,
\begin{eqnarray*}
\omega_{i,t}=\begin{cases}
 \omega_{i} &\mbox{if } t = 1, \\
\omega_{i,t-1} &\mbox{if } n-t+1 < i, \\
\frac{\|A_i^{-1}C_{i-1}\|}{1-\|A_i^{-1}B_{i}\|\omega_{i+1,t-1}} & \mbox{else}.
 \end{cases}
\end{eqnarray*}
Using these definitions we can easily prove the following modified version of Theorem~\ref{teo:blockbounds}, which refines the bounds \eqref{offdiagbound1}, \eqref{offdiagbound2} and \eqref{diagbounds} as $t$ increases, and
which generalizes~\cite[Theorems~3.4 and~3.5]{Nabben99} from the scalar to the block case.

\begin{theorem}\label{teo:iterblockbounds}
If $A$ is a matrix as in Lemma~\ref{lem:sequences} with $A^{-1}=[Z_{ij}]$, then for each $t~=~1,\ldots,n-1$,
\begin{align}
\|Z_{ij}\|&\leq \|Z_{jj}\| \prod_{k=i}^{j-1}\tau_{k,t},\quad\text{for all $i<j$},\label{offdiagbounditer1}\\
\|Z_{ij}\|&\leq \|Z_{jj}\| \prod_{k=j+1}^{i}\omega_{k,t},\quad\text{for all $i>j$}.\label{offdiagbounditer2}
\end{align}
Moreover, for $i=1,\dots,n$,
\begin{align}\label{diagboundsiter}
\frac{\|I\|}{\|A_i\|+\tau_{i-1,t}\|C_{i-1}\|+\omega_{i+1,t}\|B_{i}\|} &\leq \|Z_{ii}\| \\
&\leq \frac{\|I\|}{\|A_i^{-1}\|^{-1}-\tau_{i-1,t}\|C_{i-1}\|-\omega_{i+1,t}\|B_i\|},\nonumber
\end{align}
provided that the denominator of the upper bound is larger than zero,
and where we set $C_0=B_n=0$, and $\tau_{0,t}=\omega_{n+1,t}=0$.
\end{theorem}

Note that the statements of Theorem~\ref{teo:iterblockbounds} with $t=1$ are the same as those in Theorem~\ref{teo:blockbounds}. By construction, the sequences $\{\tau_{i,t}\}_{t=1}^{n-1}$ and $\{\omega_{i,t}\}_{t=1}^{n-1}$ are decreasing, 
and hence the bounds \eqref{offdiagbounditer1}, \eqref{offdiagbounditer2} and \eqref{diagboundsiter}
become tighter as $t$ increases. However, since we have used the submultiplicativity property of the matrix norm in the derivation, it is not guaranteed that the bounds in Theorem~\ref{teo:iterblockbounds} with $t=n-1$
will give the exact norms of the blocks of $A^{-1}$. This is a difference to the scalar case, where in the last refinement step one obtains the exact inverse; see~\cite{Nabben99}.

Finally, let us define
\begin{equation*}
\rho_{1,t}:=\max_i \tau_{i,t},\quad \omega_{2,t}:=\max_i \omega_{i,t},\quad
\mbox{for $t=1,\ldots,n-1$.}
\end{equation*}
Then the off-diagonal bounds \eqref{offdiagbounditer1} and \eqref{offdiagbounditer2} of Theorem~\ref{teo:iterblockbounds} immediately give the following result
about the decay of the norms $\|Z_{ij}\|$; cf.~\cite[Corollary~3.7]{Nabben99}

\begin{corollary}\label{cor:decay}
If $A$ is a matrix as in Theorem~\ref{teo:iterblockbounds}, then
\begin{eqnarray}
\|Z_{ij}\|\leq\rho_{1,t}^{j-i}\|Z_{jj}\|,\quad\text{ for all } i<j,\nonumber\\
\|Z_{ij}\|\leq\rho_{2,t}^{i-j}\|Z_{jj}\|,\quad\text{ for all } i>j,\nonumber
\end{eqnarray}
and for each $t~=~1,\ldots,n-1$.
\end{corollary}

\newpage
In the following we provide some numerical illustrations of the bounds in Theorem~\ref{teo:iterblockbounds}
for different values of $t$. We consider different matrices $A=[A_{ij}]$ which are row block diagonally dominant, and we compute the corresponding {matrices} $Z=[Z_{ij}]$ using the recurrences stated in Theorem~\ref{teo:Ikebe}. In all experiments we use the matrix $2$-norm, $\|\cdot\|_2$.
For each given pair $i,j$, we denote by $u_{ij}$ the value of a computed upper bound
(i.e., \eqref{offdiagbounditer1}, \eqref{offdiagbounditer2} or \eqref{diagboundsiter})
on the value $\|Z_{ij}\|_2$ and for each $i$ we denote by $l_i$ the value of the computed lower bound for the corresponding diagonal entry (i.e., \eqref{diagboundsiter}). Then the relative errors in the upper and lower bounds are given by
\begin{equation}\label{eq:RelErr}
E^{\text{u}}_{ij}=\frac{u_{ij}-\|Z_{ij}\|_2}{u_{ij}}\quad\text{ and }\quad E^{\text{l}}_{i}=\frac{\|Z_{ii}\|_2-l_{i}}{\|Z_{ii}\|_2},
\end{equation}
respectively. (Thus, both $E^{\text{u}}_{ij}$ and $E^{\text{l}}_{i}$ are between $0$ and $1$.)

\begin{figure}\begin{center}
\includegraphics[width=0.475\linewidth]{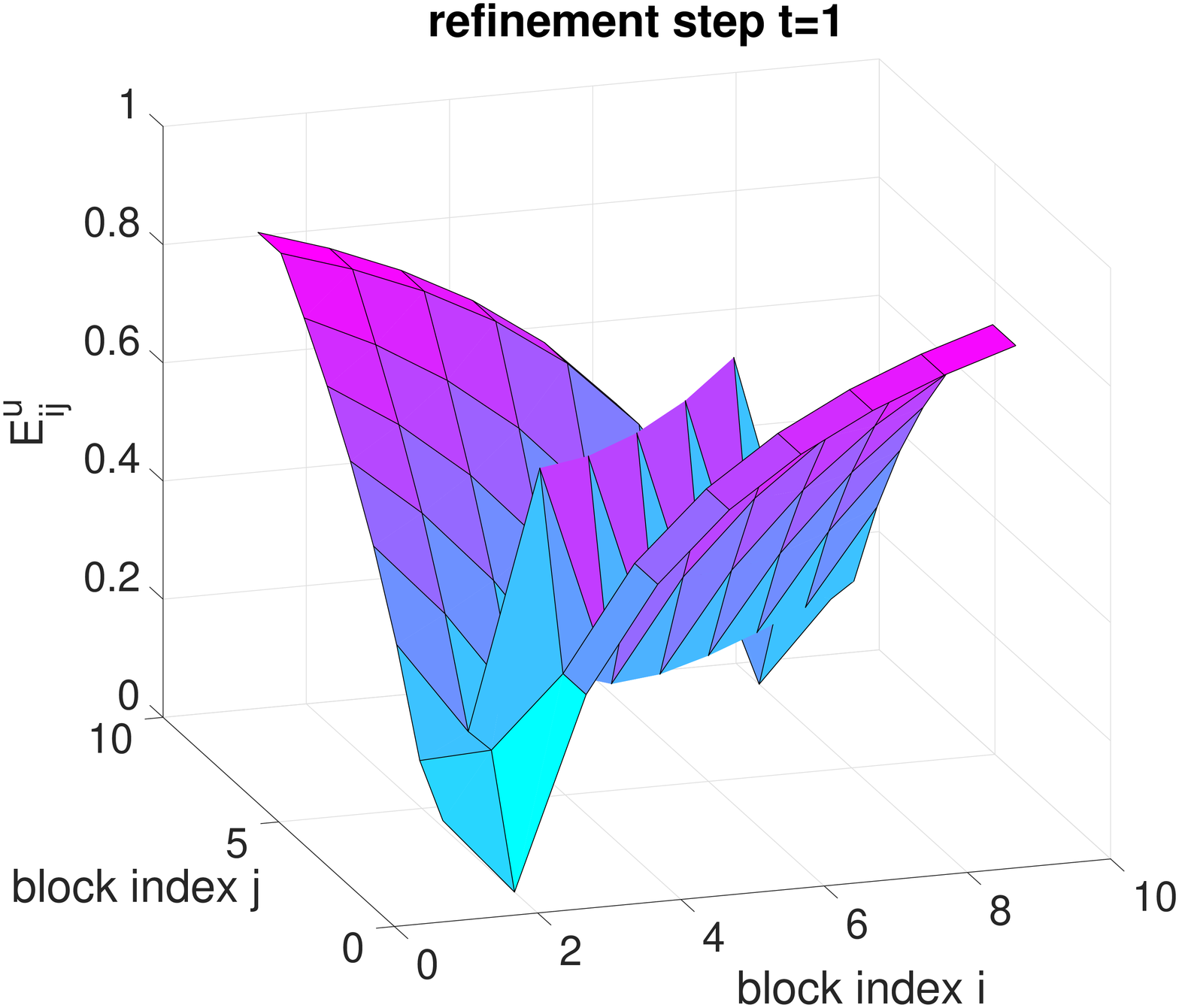}
\includegraphics[width=0.475\linewidth]{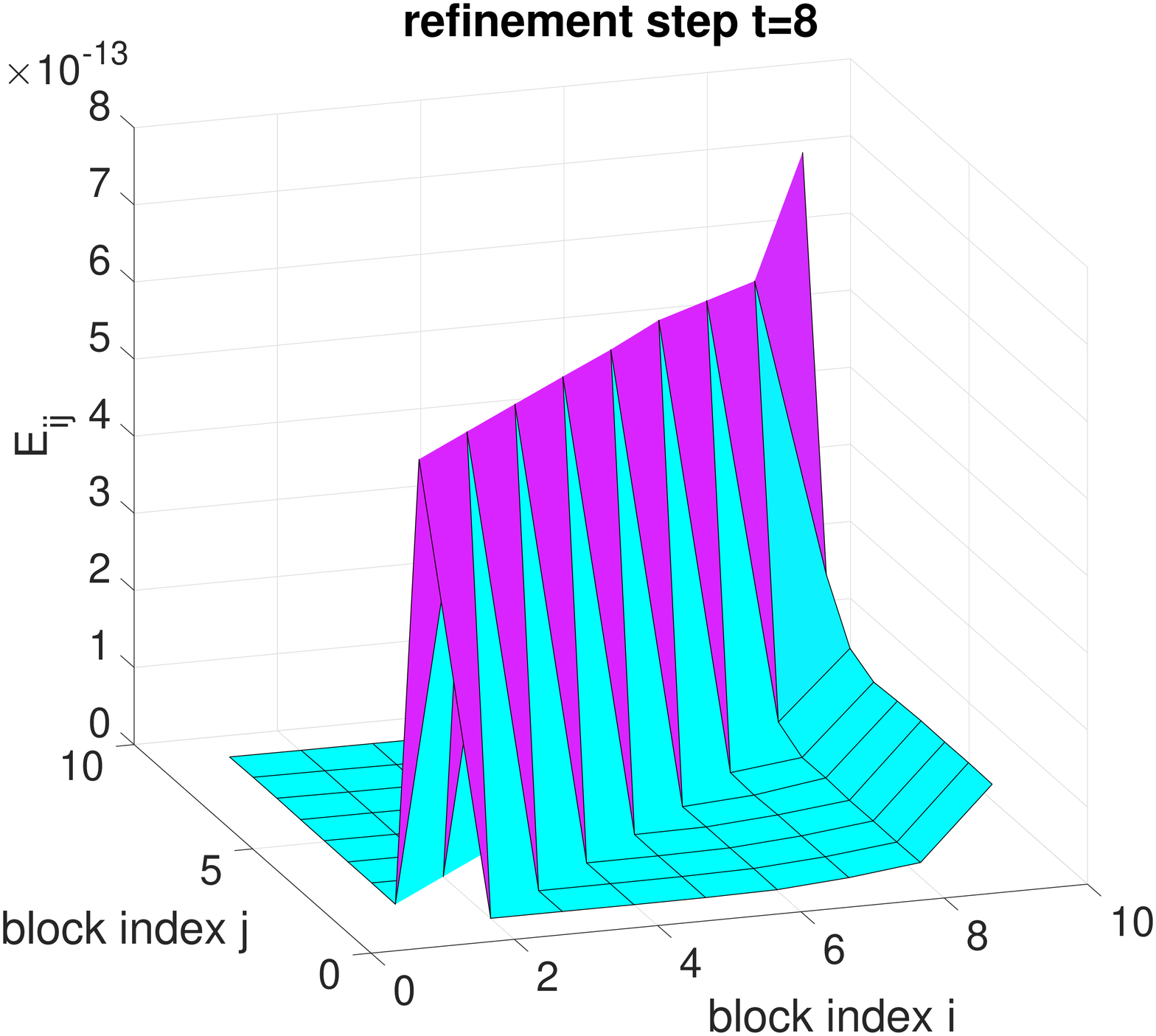}\\[1ex]
\includegraphics[width=0.475\linewidth]{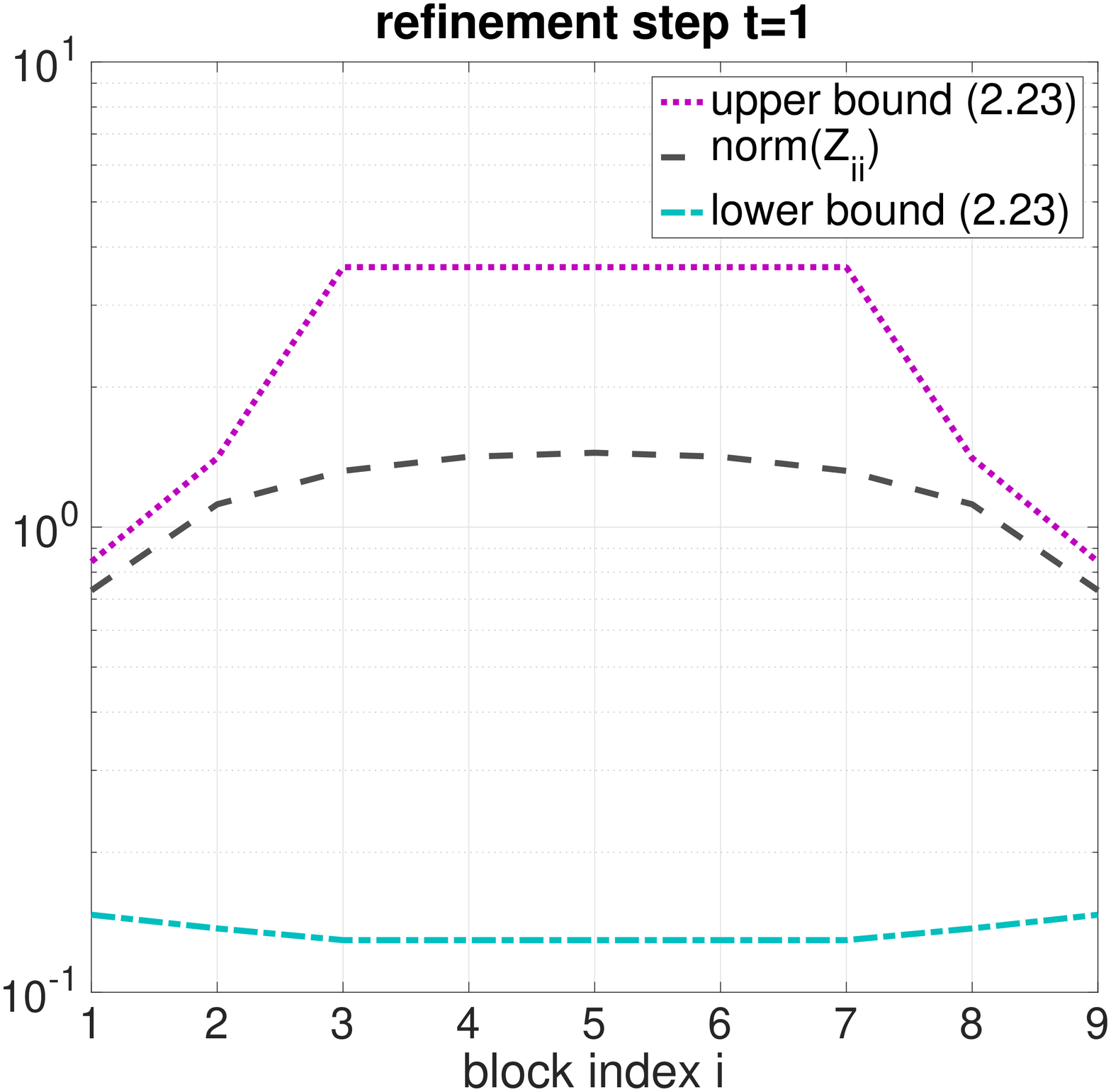}
\includegraphics[width=0.475\linewidth]{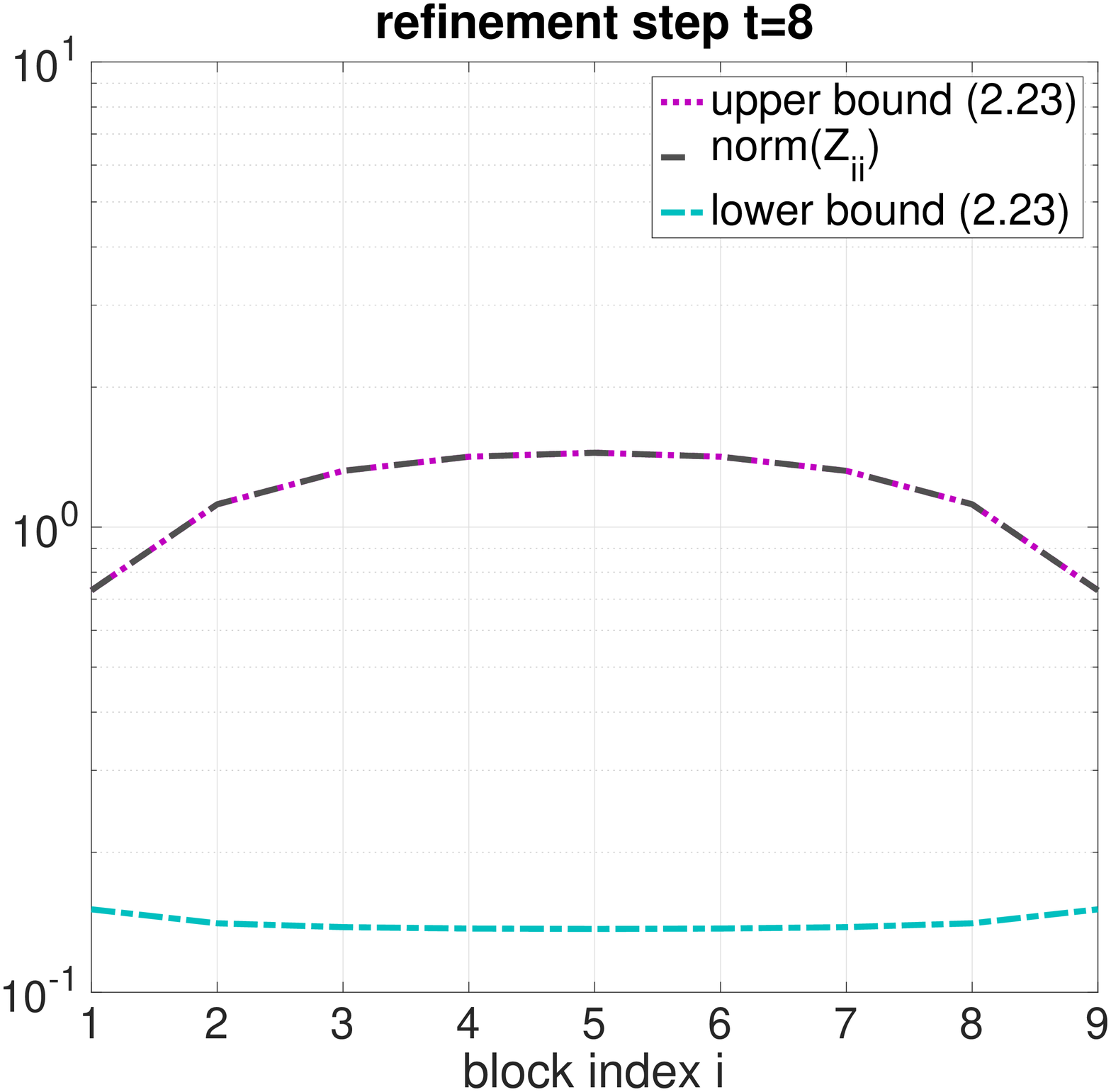}
\caption{Relative errors $E^{\text{u}}_{ij}$ (top row), upper and lower bounds on $\|Z_{ii}\|_2$ (bottom row) for the matrix $A$ of Example~\ref{ex:symm}.}\label{fig:Ex1}
\end{center}\end{figure}

\begin{example}\label{ex:symm}{\rm
We start with the symmetric block Toeplitz matrix
\begin{equation}\label{eq:kronA}
A~=T\otimes I+ I\otimes~T \;\;\in\;\;\mathbb{R}^{81\times 81},
\end{equation}
where $T=\text{tridiag}(-1,2,-1)\in\mathbb{R}^{9\times 9}$, i.e., $A$ is of the form \eqref{eq:blocktridiag} with \linebreak
$A_i=\text{tridiag}(-1,4,-1)$, and $B_i=C_i=\text{diag}(-1)$ for all $i$.
We have \linebreak $\kappa_2(A_1)~=~58.4787$, i.e., the matrix $A_1$ is quite well conditioned. For the
computed matrix $Z=[Z_{ij}]$ we obtain $\|ZA-I\|_2=2.7963\times 10^{-10}$, suggesting that
$Z$ is a reasonably accurate approximation of the exact inverse $A^{-1}$.

In the top row of Figure~\ref{fig:Ex1} we show the relative errors $E^{\text{u}}_{ij}$
for the refinement step $t=1$ (no refinement) and $t=8$ (maximal refinement). We observe that the upper bounds are quite tight already for $t=1$, and that for $t=8$ the maximal relative error is on the order $10^{-13}$, i.e., the value of the upper bound is almost exact.
In the bottom row of Figure~\ref{fig:Ex1} we show the values $\|Z_{ii}\|_2$ for
$i=1,\dots,9$, and the corresponding upper and lower bounds \eqref{diagboundsiter} for the
refinement steps $t=1$ and $t=8$. We observe that while the upper bounds on $\|Z_{ii}\|_2$
for $t=8$ almost exactly match the exact values, the lower bounds do not
improve by the iterative refinement. The maximal error of the lower bounds for the diagonal block
entries of $Z$ in the maximal refinement step is on the order $10^{-1}$.
The maximal relative errors in the upper and lower bounds and all refinement steps are shown in the following table:
%
\begin{table}[h!]\scriptsize\centering%
\begin{tabular}{c c c c c c c c c}
\hline
 t & 1 & 2 & 3 & 4 & 5 & 6 & 7 & 8 \\
 $\max_{ij}E^{u}_{ij}$ &  $0.84478$ &   $0.63381$  &  $0.39537$ &  $0.20899$  & $0.09596$ &  $0.03780$ & $0.01109$ & $7.141\times 10^{-13}$\\
 $\max_i E^{l}_{i}$ & $0.91039$ & $0.90877$ &  $0.90765$ & $0.90529$ & $0.90529$ & $0.90529$ & $0.90529$ &  $0.90529$\\
\hline
\end{tabular}
\end{table}
}\end{example}

\begin{figure}\begin{center}
\includegraphics[width=0.475\linewidth]{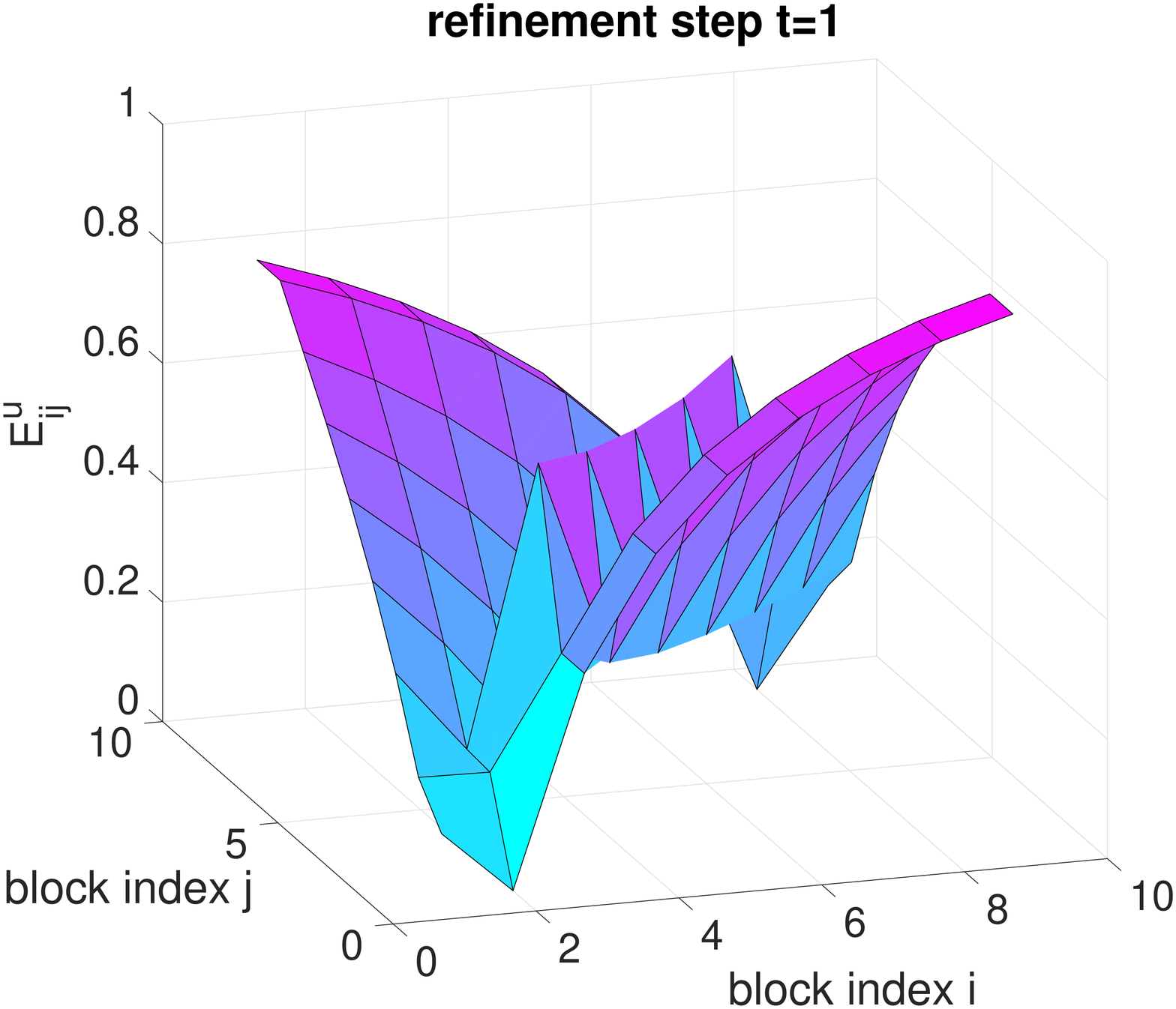}
\includegraphics[width=0.475\linewidth]{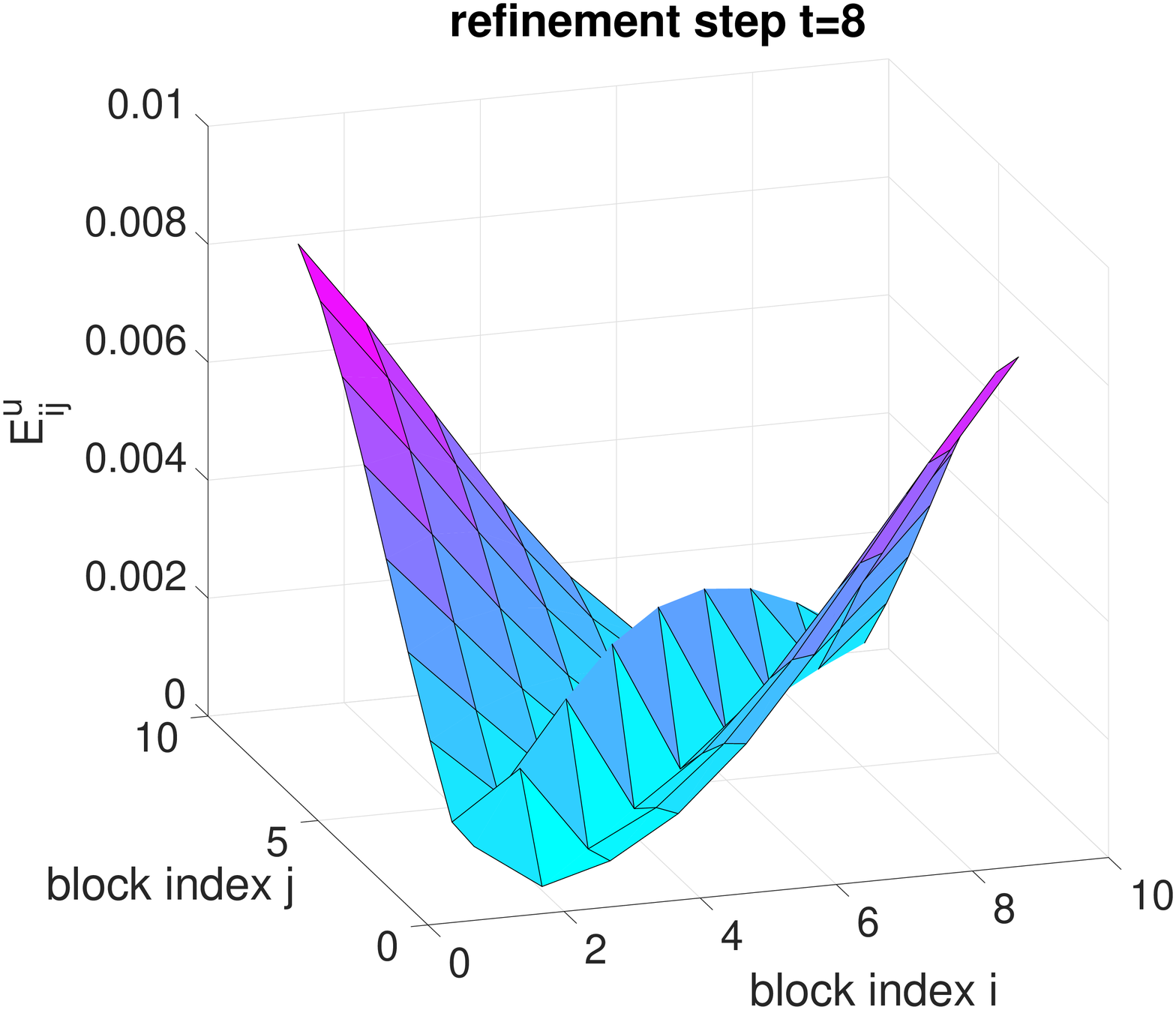}\\[1ex]
\includegraphics[width=0.475\linewidth]{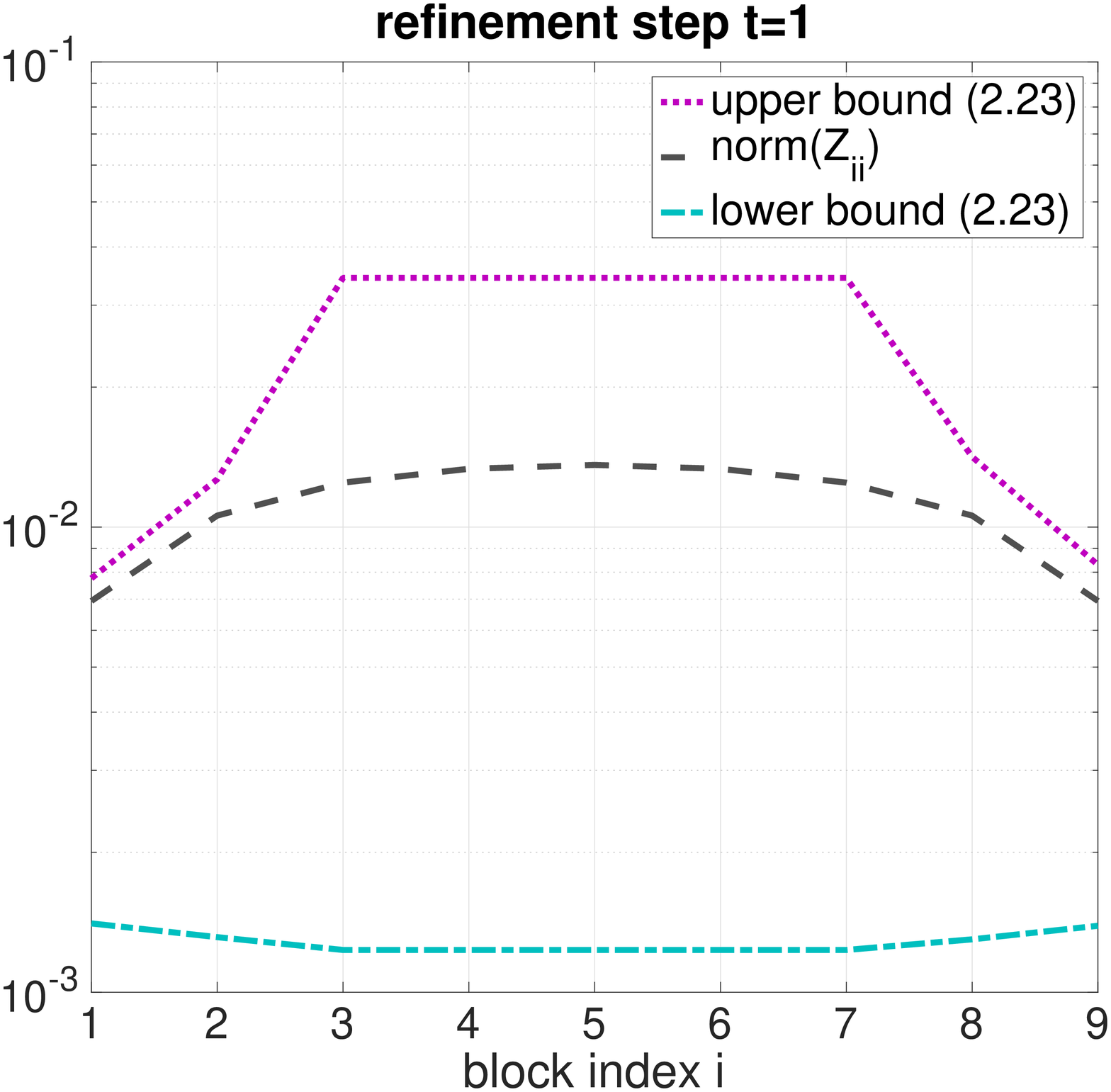}
\includegraphics[width=0.475\linewidth]{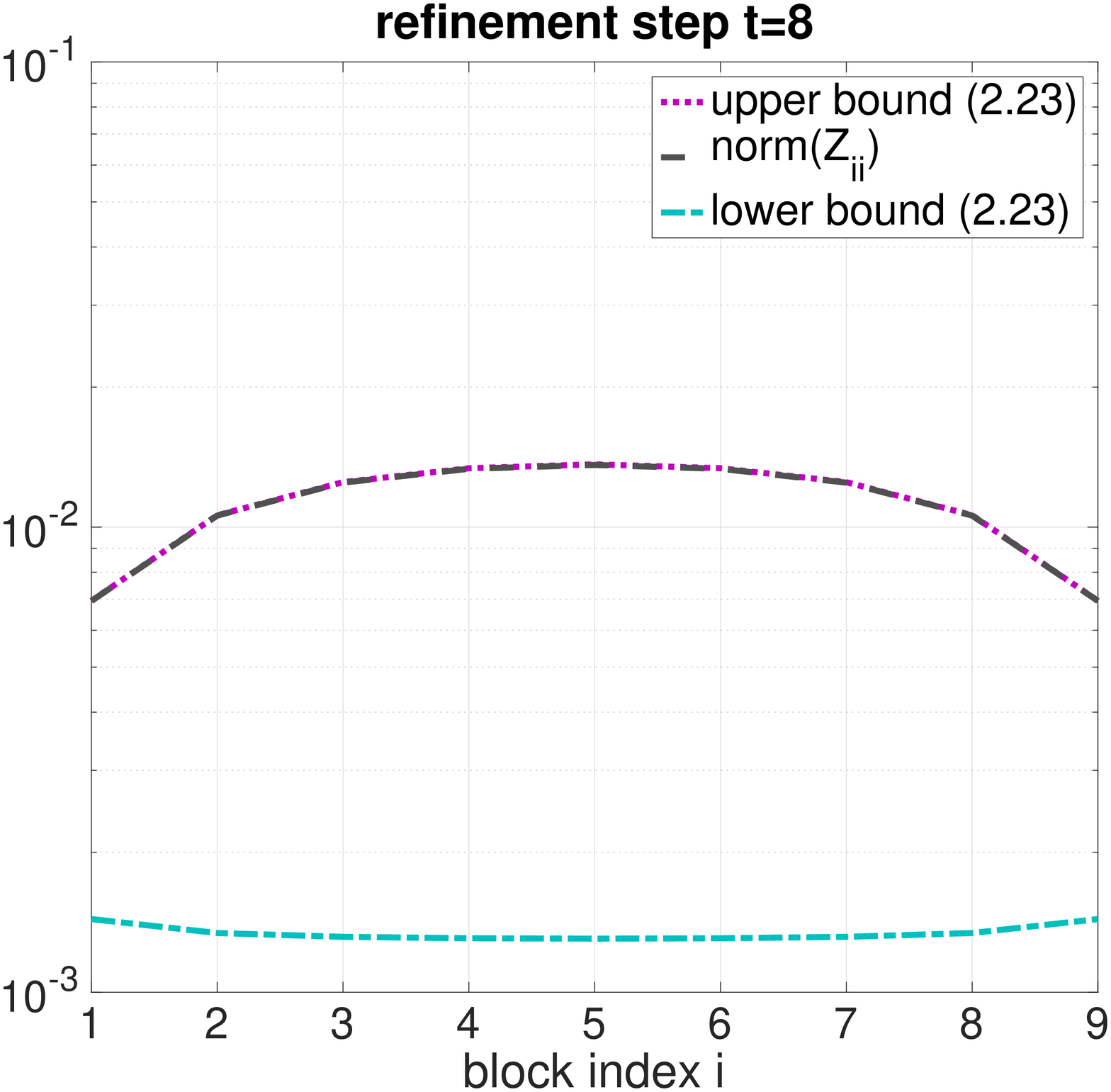}
\caption{Relative errors $E^{\text{u}}_{ij}$ (top row), and upper and lower bounds on $\|Z_{ii}\|_2$
(bottom row) for the matrix $A$ of Example~\ref{ex:nonsym1}.}\label{fig:Ex2}
\end{center}\end{figure}

\newpage
\begin{example}\label{ex:nonsym1}{\rm
Let $A$ be the nonsymmetric block Toeplitz matrix of the form \eqref{eq:kronA} with
$T=\text{tridiag}(-110,209.999,-99.999)\in\mathbb{R}^{9\times 9}$, i.e., $A$
again takes the form \eqref{eq:blocktridiag} with $A_i=\text{tridiag}(-110,419.999,-99.999)$, $B_i=\text{diag}(-110)$, and $C_i=\text{diag}(-99.999)$. The condition number in this case is $\kappa_2(A)=57.5725$,  and for the computed matrix $Z$ we obtain $\|ZA-I\|_2=1.5151\times 10^{-10}$.

The top row of Figure~\ref{fig:Ex2} shows the relative errors for the refinement steps $t=1$ and $t=8$.
We observe that for this nonsymmetric example the upper bounds are not as accurate as those given in the symmetric case, producing a maximal relative error at refinement step $t=8$ on the order $10^{-3}$.
The bottom row of Figure~\ref{fig:Ex2} shows the upper and lower bounds \eqref{diagboundsiter}
as well as the values $\|Z_{ii}\|_2$ for $i=1,\dots,9$, and refinement steps $t=1$ and $t=8$.
Again we can observe that while we obtain a reasonable approximation in the upper bounds on $\|Z_{ii}\|_2$
for $t=8$, the lower bounds almost do not improve by the iterative refinement process.
The maximal relative errors in the upper and lower bounds and all refinement steps is shown in the
following table:
%
%
\begin{table}[h!]\scriptsize\centering
\begin{tabular}{c c c c c c c c c}
\hline
  t & 1 & 2 & 3 & 4 & 5 & 6 & 7 & 8 \\
$\max_{ij}E^{u}_{ij}$ &  $0.88856$ & $ 0.70640$ & $0.46700$ &  $0.25859$ & $0.12442$ & $0.05378$ &  $0.02140$ &  $0.00824$\\
$\max_i E^{l}_{i}$ & $0.90934$ & $0.90768$ & $0.90652$ & $0.90411$ & $0.90411$ & $0.90411$ &  $0.90411$ &  $0.90411$\\
\hline
\end{tabular}
\end{table}

}\end{example}

\begin{figure}\begin{center}
\includegraphics[width=0.475\linewidth]{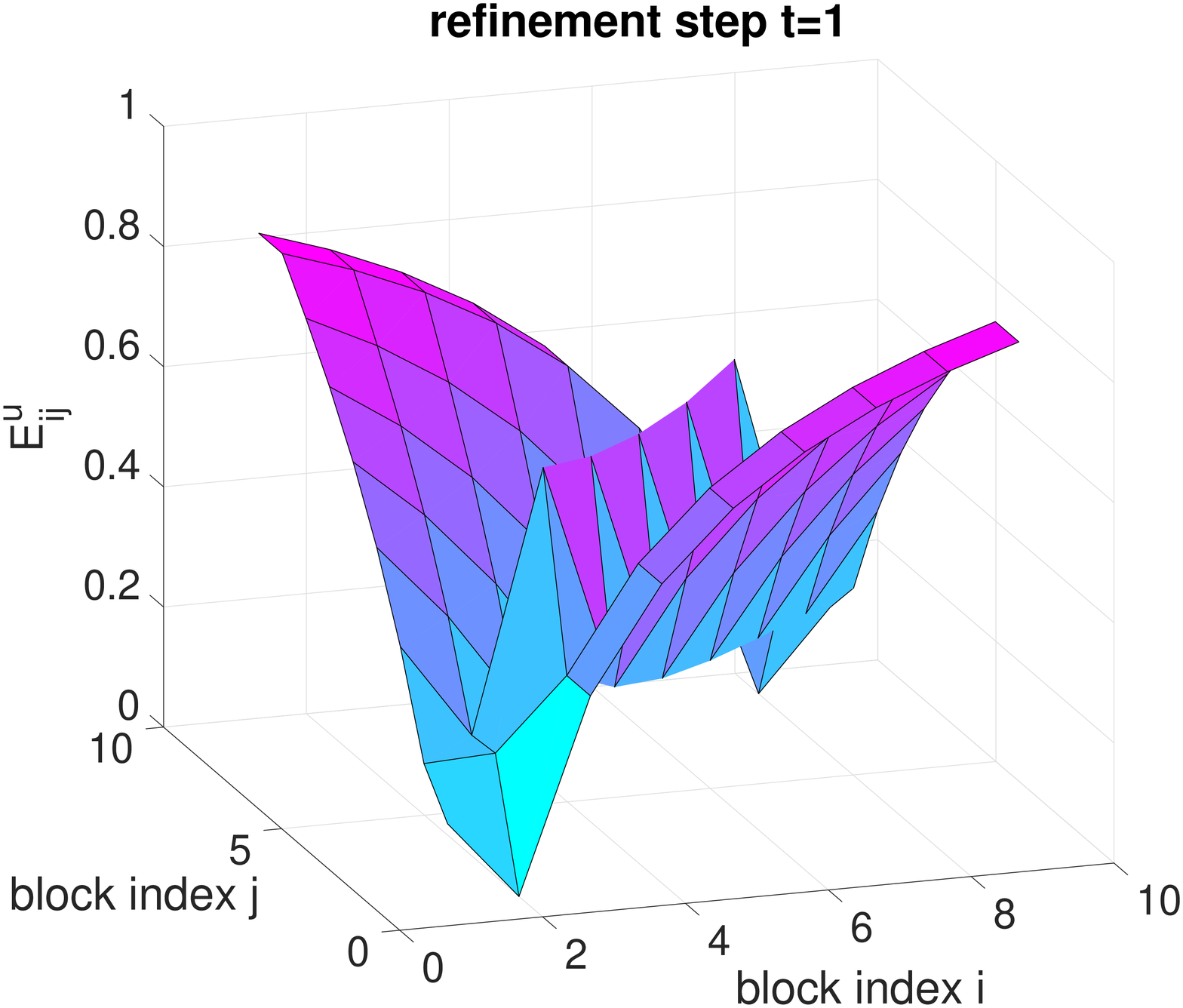}
\includegraphics[width=0.475\linewidth]{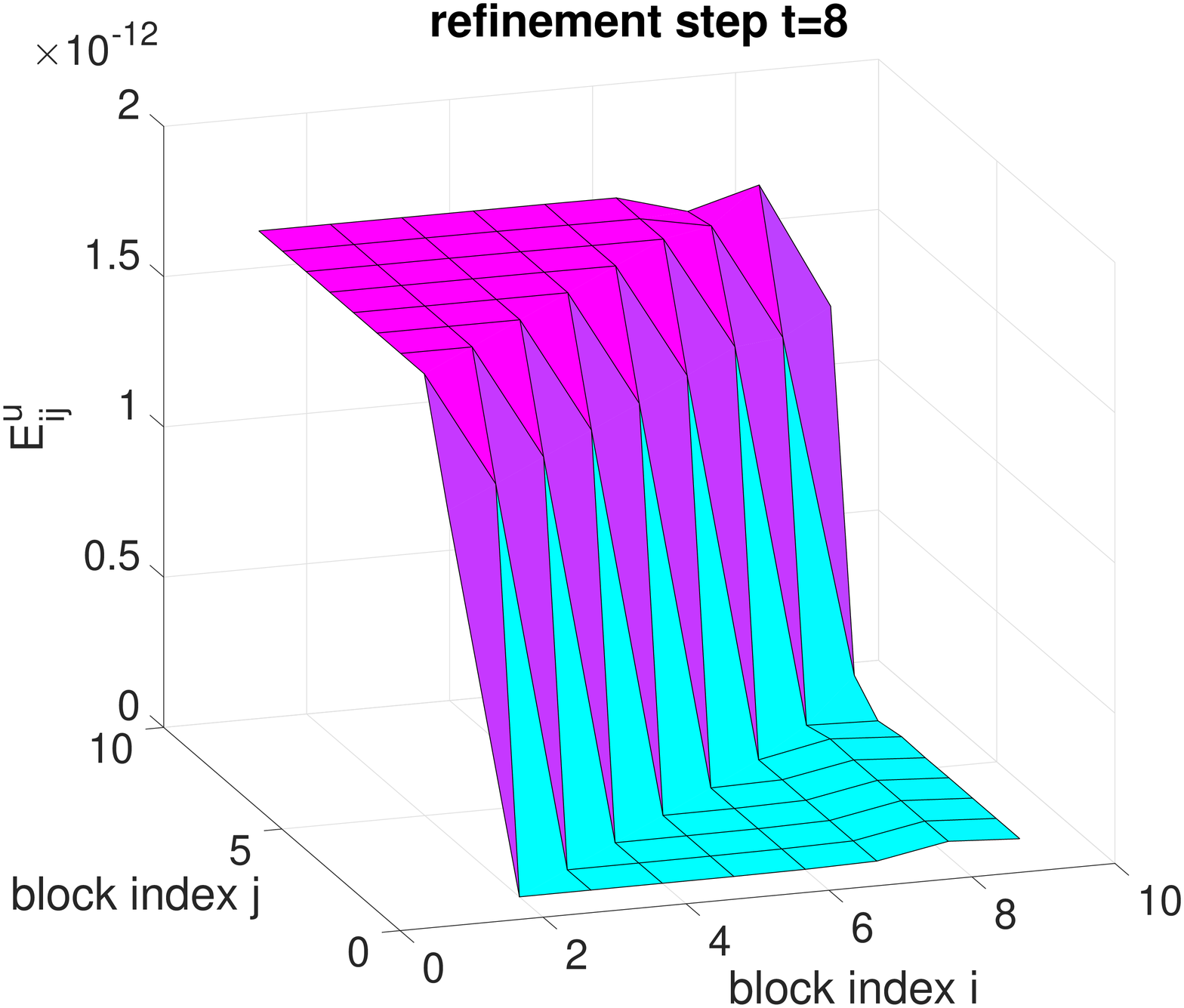}\\[1ex]
\includegraphics[width=0.475\linewidth]{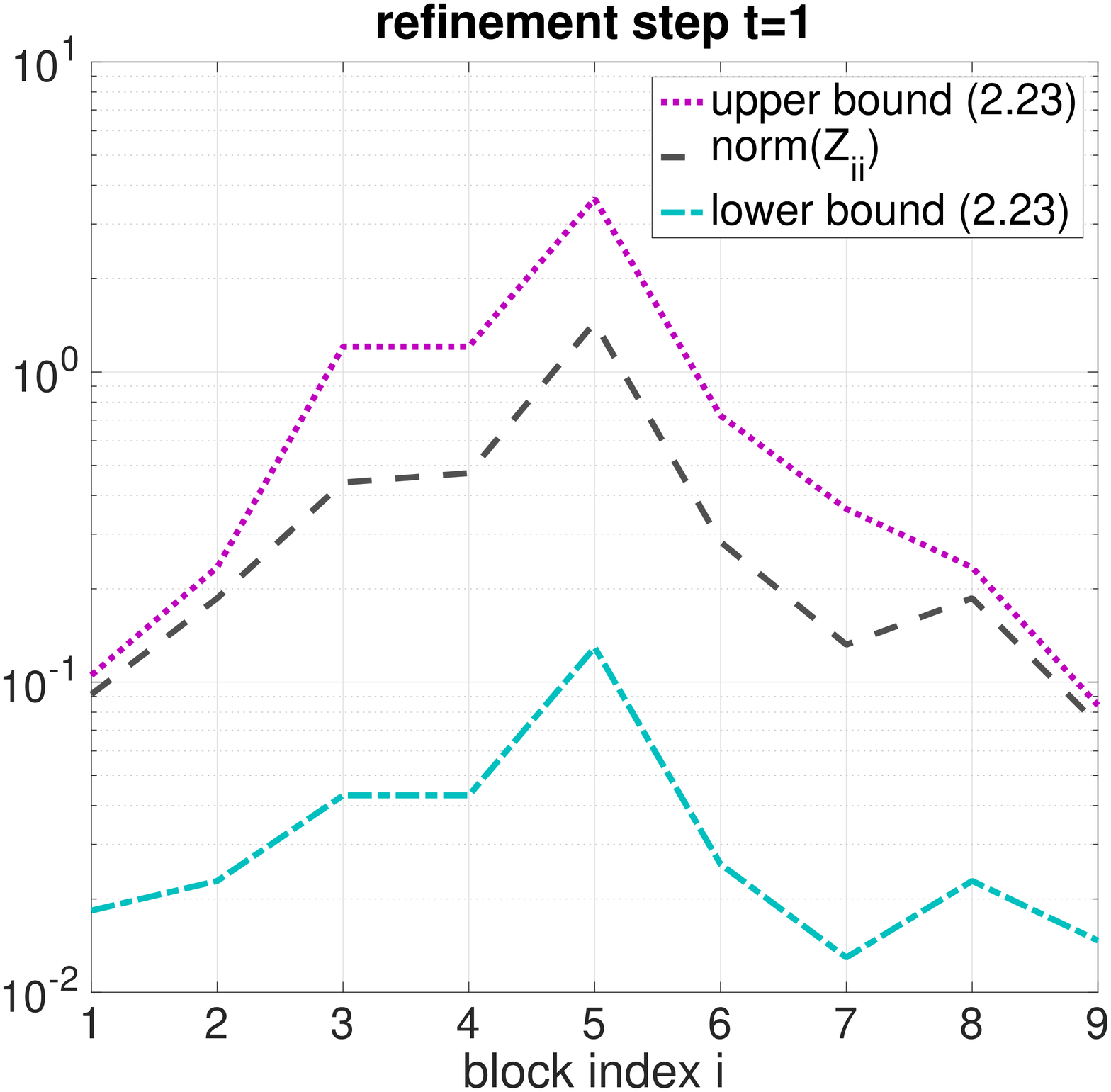}
\includegraphics[width=0.475\linewidth]{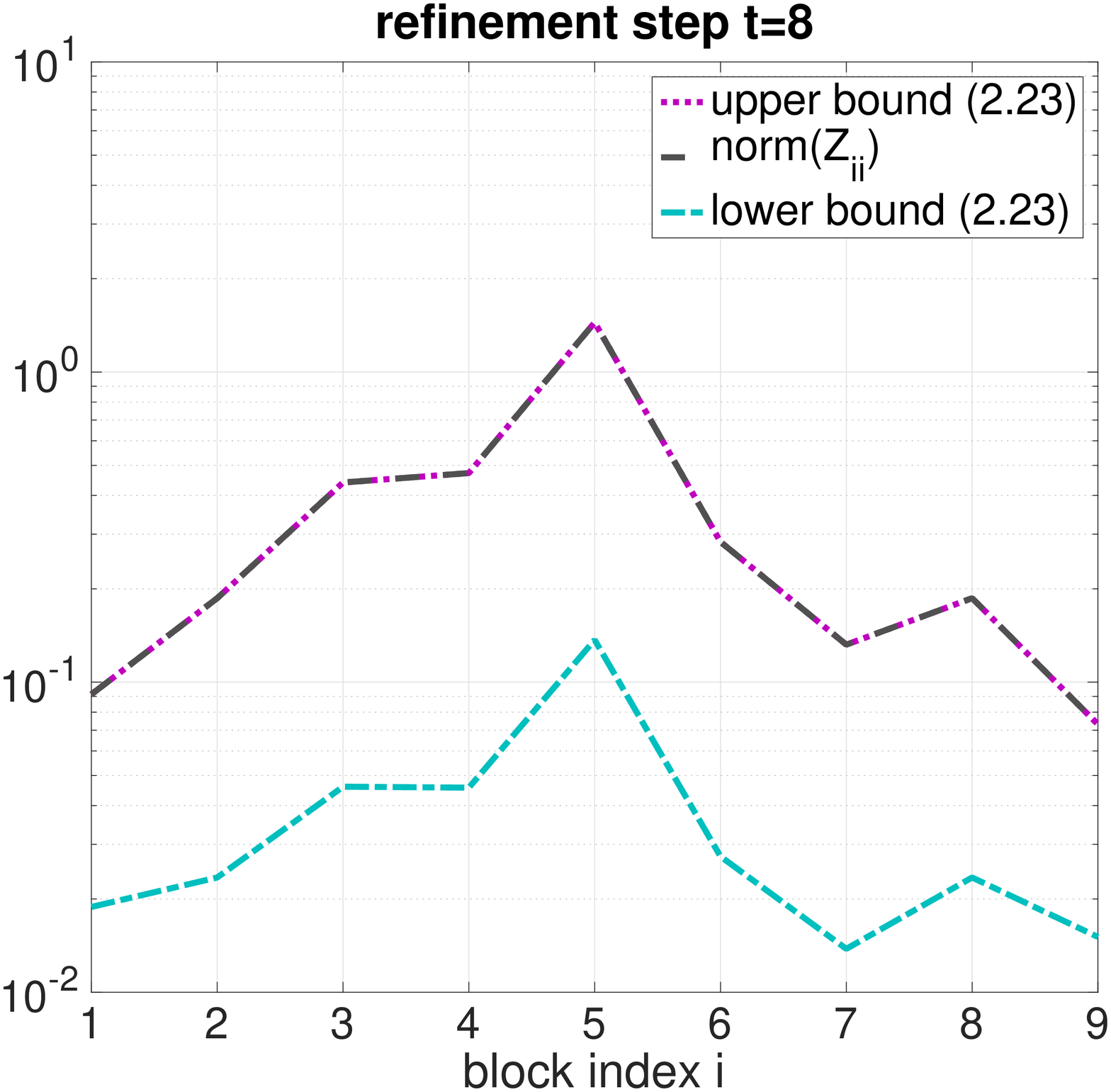}
\caption{Relative errors $E^{\text{u}}_{ij}$ (top row), and upper and lower bounds on $\|Z_{ii}\|_2$
(bottom row) for the matrix $A$ of Example~\ref{ex:nonsym2}.}\label{fig:Ex3}
\end{center}\end{figure}

\pagebreak
\begin{example}\label{ex:nonsym2}{\rm
We now consider the  nonsymmetric block tridiagonal matrix
\[A~=(R\otimes I)(T\otimes I+ I\otimes~T) \;\;\in\;\;\mathbb{R}^{81\times 81},\]
where $T$ is given as in Example \ref{ex:symm}, and $R\in\mathbb{R}^{9\times 9}$ is a random
diagonal matrix with nonzero integer entries between $0$ and $10$
and constructed in MATLAB with the command \texttt{R = diag(ceil(10*rand(9,1)))}.
Thus, $A$ is of the form \eqref{eq:blocktridiag} with random tridiagonal Toeplitz matrices
$A_i$, and random constant diagonal matrices $B_i$ and $C_i$ for all $i$.
For this matrix we have $\kappa_2(A)=518.9988$, and the computed matrix $Z$ yields
$\|ZA-I\|_2=1.0519\times 10^{-9}$. The relative errors in the bounds are shown in
Figure~\ref{fig:Ex3} and in the following table:
\begin{table}[h!]\scriptsize\centering
\begin{tabular}{c c c c c c c c c}
\hline
 t & 1 & 2 & 3 & 4 & 5 & 6 & 7 & 8 \\
 $\max_{ij}E^{u}_{ij}$ &  $0.84477$ & $0.63381$  &  $0.39537$ &  $0.20898$ & $0.09595$ &  $0.03780$ &  $0.01109$ &  $1.313\times 10^{-12}$\\
 $\max_{i}E^{l}_{i}$ &  $0.91039$ & $0.90877$ & $0.90765$ & $0.90529$ & $0.90529$ &  $0.90529$ &  $0.90529$ &  $0.90529$\\
  \hline
\end{tabular}
\end{table}
}\end{example}

\newpage

\begin{figure}\begin{center}
\includegraphics[width=0.475\linewidth]{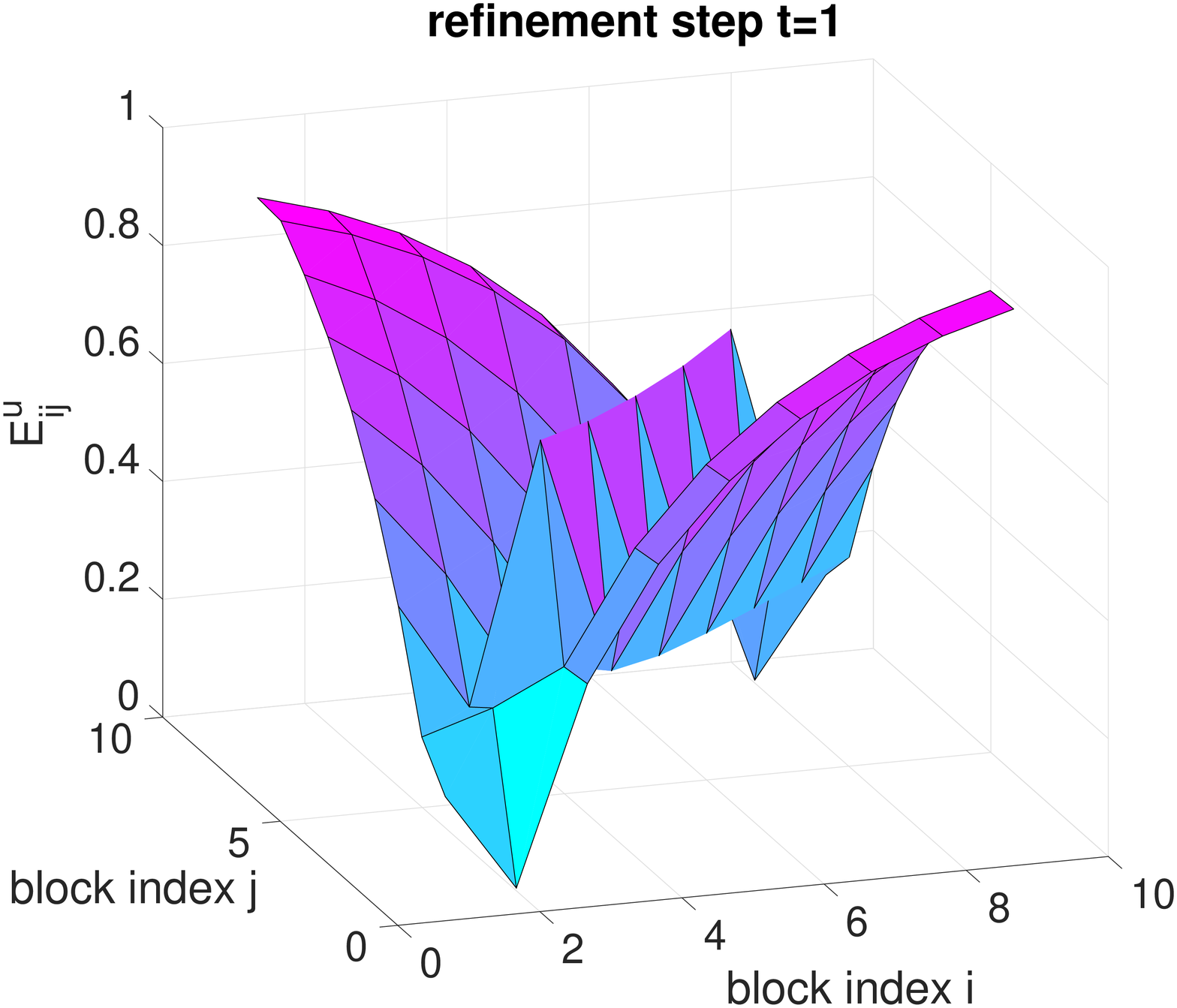}
\includegraphics[width=0.475\linewidth]{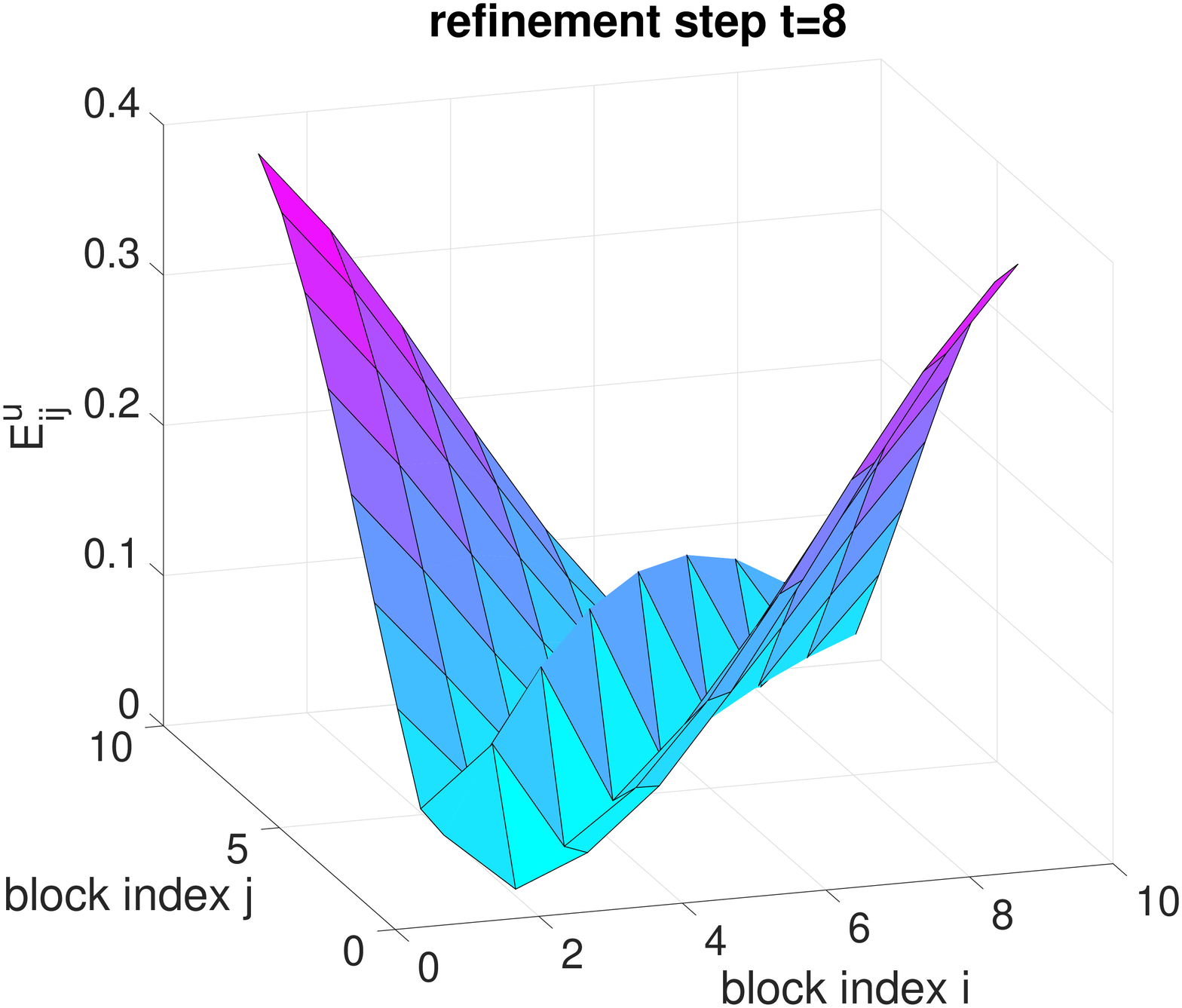}\\[1ex]
\includegraphics[width=0.475\linewidth]{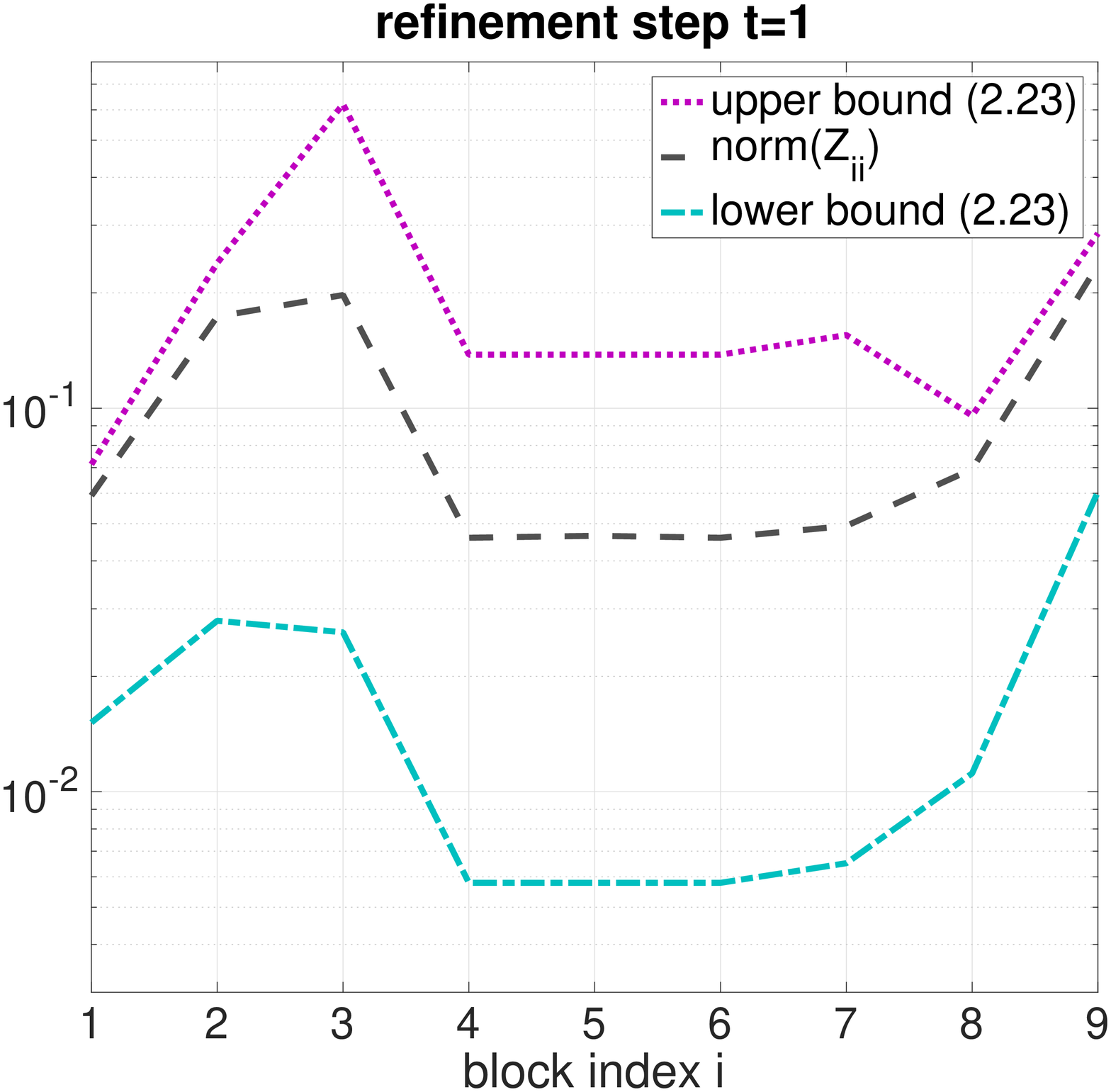}
\includegraphics[width=0.475\linewidth]{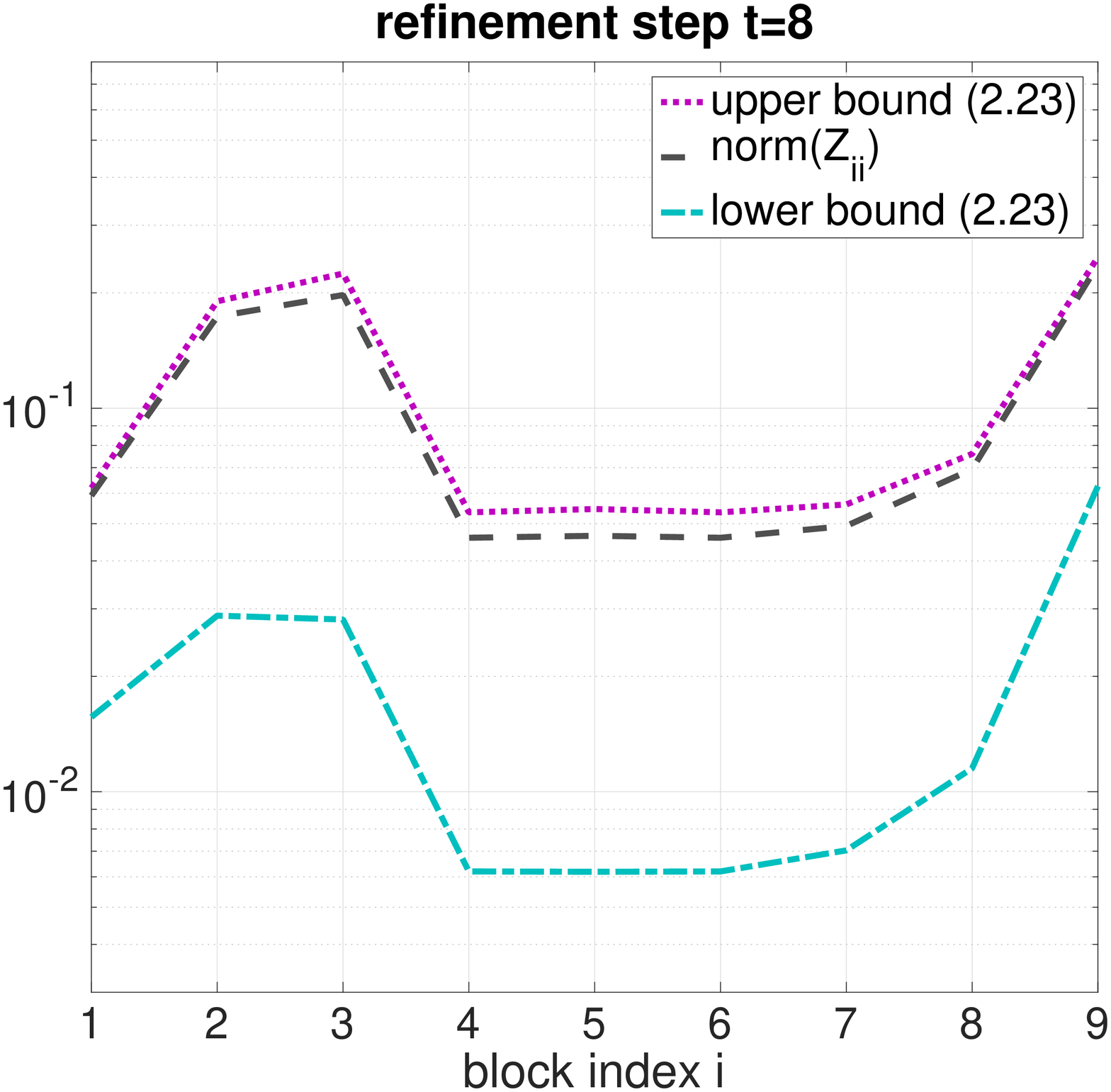}
\caption{Relative errors $E^{\text{u}}_{ij}$ (top row), and upper and lower bounds on $\|Z_{ii}\|_2$
(bottom row) for the matrix $A$ of Example~\ref{ex:nonsym3}.}\label{fig:Ex4}
\end{center}\end{figure}

\begin{example}\label{ex:nonsym3}{\rm
Finally, we consider the nonsymmetric block tridiagonal matrix
\begin{equation*}
A=\\(R\otimes I)\tridiag(\tridiag(-0.01,-2,1),\tridiag(-2,10,-2),\tridiag(-0.01,-2,1)),
\end{equation*}
with $A\in\mathbb{R}^{81\times 81}$, and where $R\in\mathbb{R}^{9\times 9}$ is a random diagonal matrix constructed as in Example \ref{ex:nonsym2}. In this case $A$ takes the form \eqref{eq:blocktridiag} with $A_i$, $B_i$ and $C_i$ random tridiagonal Toeplitz matrices with integer entries for all $i$.
For this matrix we have $\kappa_2(A)=162.376$, and $\|ZA-I\|_2=3.7668\times 10^{-9}$.
The relative errors in the bounds are shown in Figure~\ref{fig:Ex4} and the following table:
\begin{table}[h!]\scriptsize\centering
\begin{tabular}{c c c c c c c c c}
\hline
 t & 1 & 2 & 3 & 4 & 5 & 6 & 7 & 8 \\
 $\max_{ij}E^{u}_{ij}$ &  $0.90505$ & $0.7769$ & $0.63219$ & $0.51856$ & $0.44929$ & $0.41345$ &  $0.39690$ & $0.38998$\\
 $\max_{i}E^{l}_{i}$ &  $0.87551$ & $0.87284$ & $0.87099$ & $0.86698$ & $0.86698$ & $0.86698$ &  $0.86698$ & $0.86698$\\
  \hline
\end{tabular}
\end{table}

}\end{example}

\newpage
\section{Inclusion regions for eigenvalues}\label{s:eigenregion}
In this section we generalize a result of Feingold and Varga on eigenvalue inclusion regions of block matrices. We start with the following generalization
of~\cite[Theorem~1]{FeiVar62}; also cf.~\cite[Theorem~6.2]{Var2011}.




\begin{lemma}\label{lem:nonsingular}
If a matrix $A$ as in \eqref{eq:blockmat} is row block strictly diagonally dominant,
then $A$ is nonsingular.
\end{lemma}

\emph{Proof.}
The proof closely follows the proof of~\cite[Theorem~1]{FeiVar62}. Suppose that $A$ is row
block strictly diagonally dominant but singular. Then there exists a nonzero block vector $X$,
partitioned conformally with respect to the partition of~$A$ in~\eqref{eq:blockmat}, such that
\[
A\left[\begin{array}{c}X_1\\\vdots\\X_n\end{array}\right]=0.
\]
This is equivalent to
\[
A_{ii}X_i+\sum_{\atopfrac{j=1}{j\neq i}}^{n} A_{ij}X_j=0,\quad 1\leq i \leq n,
\]
and, since the diagonal blocks $A_{ii}$ are nonsingular,
\[
X_i=-\sum_{\atopfrac{j=1}{j\neq i}}^{n} A_{ii}^{-1}A_{ij}X_j,\quad 1\leq i \leq n,
\]
Without loss of generality we can assume that $X$ is normalized such that $\|X_i\|~\leq~1$
for all $1\leq i \leq n$, with equality for some $i=r$. For this index we obtain
\[
1=\|X_r\|=\Big\|\sum_{\atopfrac{j=1}{j\neq i}}^{n} A_{rr}^{-1}A_{rj}X_j\Big\|\leq \sum_{\atopfrac{j=1}{j\neq r}}^{n}
\|A_{rr}^{-1}A_{rj}\|\|X_j\|\leq \sum_{\atopfrac{j=1}{j\neq r}}^{n} \|A_{rr}^{-1}A_{rj}\|,
\]
which contradicts the assumption that $A$ is row block strictly diagonally dominant. Thus, $A$
must be nonsingular. \qed

If $\lambda$ is an eigenvalue of $A$, then $A-\lambda I$ is singular, and hence $A-\lambda I$ cannot be block strictly diagonally dominant. This immediately gives the following result, which generalizes~\cite[Theorem~2]{FeiVar62}; also cf.~\cite[Theorem~6.3]{Var2011}.

\begin{corollary}\label{cor:inclusion}
If a matrix $A$ is as in \eqref{eq:blockmat}, and $\lambda$ is an eigenvalue of $A$, then
there exists at least one $i\in\{1,\dots,n\}$ with
\begin{equation}\label{eq:evals}
\sum_{\atopfrac{j=1}{j\neq i}}^{n} \| (A_{ii}-\lambda I)^{-1} A_{ij} \|\geq 1.
\end{equation}
\end{corollary}

If all the blocks $A_{ij}$ of $A$  are of size $1\times 1$, and $\|A_{ij}\|=|A_{ij}|$, then this result reduces to the classical Gershgorin Circle Theorem.

Corollary \ref{cor:inclusion} shows that each eigenvalue $\lambda$ of $A$ must be contained in the union of the sets
\[
G_i^{new}=\Big\{z\in\mathbb{C} : \sum_{\atopfrac{j=1}{j\neq i}}^{n} \| (A_{ii}-z I)^{-1} A_{ij}\|\geq 1\Big\},
\]
for $i=1,\dots,n$. Due to the submultiplicativity property of the matrix norm, the sets $G_i^{new}$ are potentially smaller than the ones proposed in~\cite[Definition~3]{FeiVar62},
\begin{small}\[
G_i^{FV}=\Big\{z\in\mathbb{C} : \sum_{\atopfrac{j=1}{j\neq i}}^{n} \| (A_{ii}-z I)^{-1}\| \|A_{ij}\|\geq 1\Big\},
\]\end{small}
i.e., we have $G_i^{new}\subseteq G_i^{FV}$. We will illustrate this fact with numerical examples.

\begin{figure}\begin{center}
\includegraphics[width=0.49\linewidth]{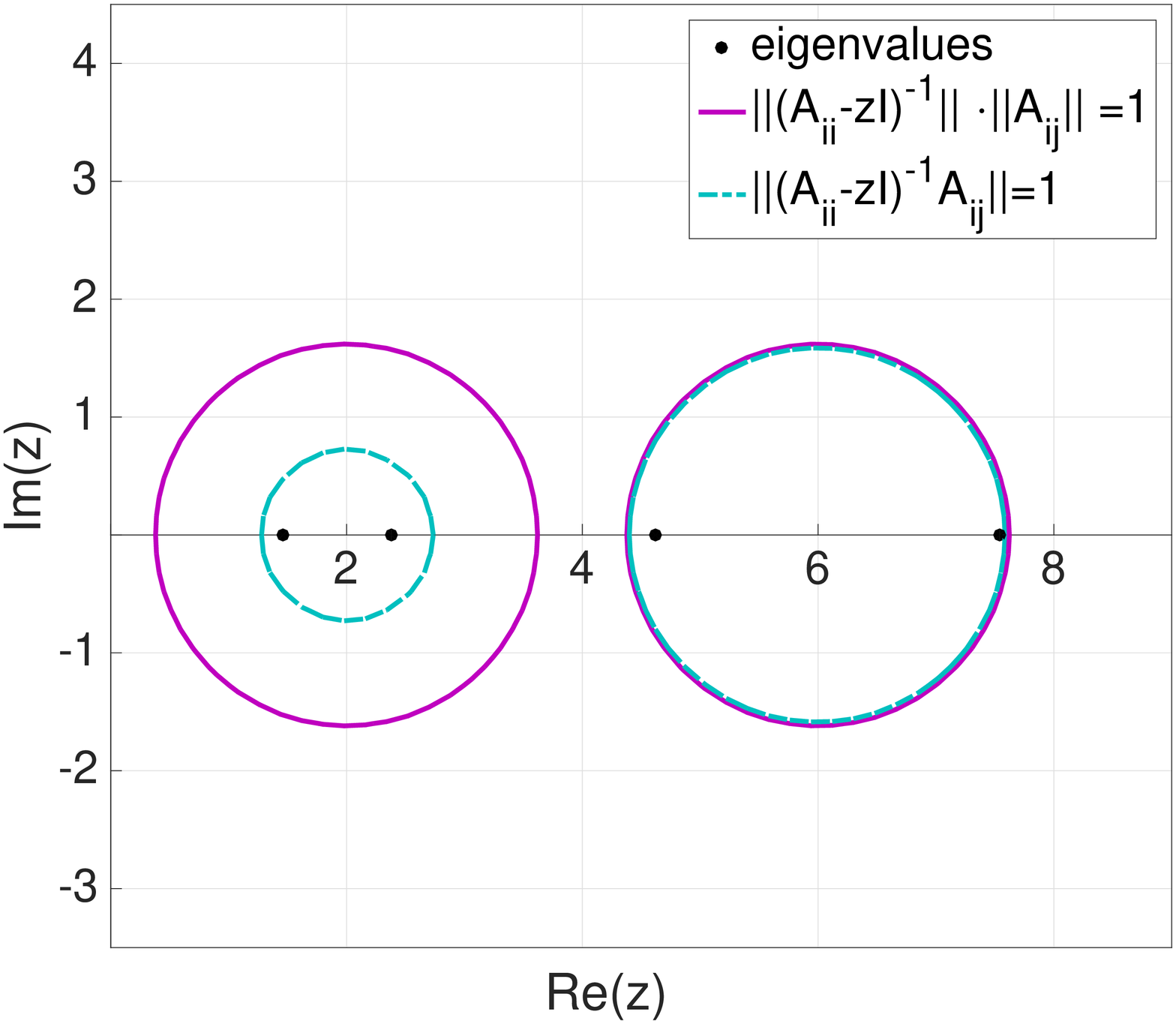}
\includegraphics[width=0.49\linewidth]{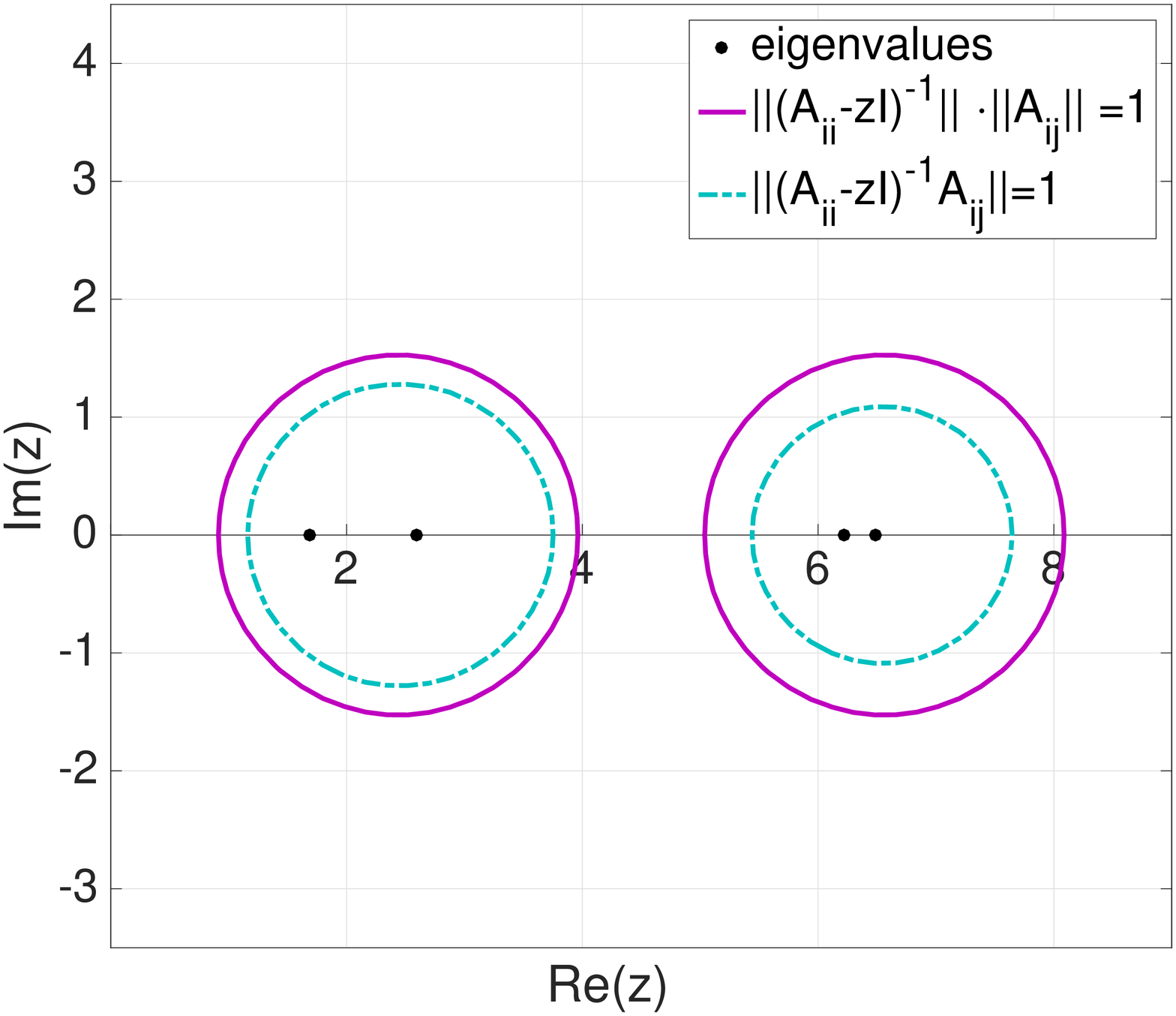}
\caption{Eigenvalue inclusion regions obtained from the sets
$G_i^{new}$ and $G_i^{FV}$ in Example~\ref{ex:inclusion}.}\label{fig:cont}
\end{center}\end{figure}

\begin{example}\label{ex:inclusion}{\rm
We first consider the symmetric matrix
\begin{equation*}
A=
\left[ \begin{array}{rr|rr}
4	& -2	& -1 & 1	\\
-2	&	4	&  0 & -1	\\ \hline
-1	& 0	&	4	&-2	\\
1	& -1	& -2	& 4
\end{array} \right]
=
\left[ \begin{array}{c|c}
A_{11}	& A_{12}	\\ \hline
A_{21}	& A_{22}
\end{array} \right],
\end{equation*}
which has the eigenvalues $1.4586$, $2.3820$, $4.6180$, and $7.5414$ (computed in MATLAB
and rounded to five significant digits). The left part of Figure~\ref{fig:cont} shows the boundaries of the corresponding sets $G_j^{new}$ and $G_j^{FV}$ for $j=1,2$, i.e., the curves for $z\in {\mathbb C}$ where
%

\begin{equation*}
\| (A_{ii}-z I)^{-1} A_{ij} \|=1,\text{ and }
\| (A_{ii}-z I)^{-1}\|\| A_{ij} \|=1,
\quad i,j\in\{1,2\},\quad i\neq j,
\end{equation*}
respectively. Clearly, the sets $G_i^{new}$ give tighter inclusion regions for
the eigenvalues than the sets $G_i^{FV}$ as well as the usual Gershgorin circles for the matrix $A$, which are given by the two circles centered at $z=4$ of radius~3 and~4.

We next consider the nonsymmetric matrix
\begin{equation}\label{eq:ELNmatrix2}
A=
\left[ \begin{array}{rr|rr}
4	& -2	& -0.5 & 0.5	\\
-2	&	5	&  -1.4 & -0.5	\\ \hline
-0.5	& 0	&	4	&-2	\\
0.5	& -0.5	& -2	& 4
\end{array} \right]
=
\left[ \begin{array}{c|c}
A_{11}	& A_{12}	\\ \hline
A_{21}	& A_{22}
\end{array} \right],
\end{equation}
which has the eigenvalues $1.6851$, $2.5959$, $6.2263$, and $6.4927$. As shown in the right part of Figure~\ref{fig:cont}, the sets $G_i^{new}$ again give tighter inclusion regions than the sets $G_i^{FV}$ as well as the usual Gershgorin circles.
}\end{example}



\section*{Acknowledgments}
Carlos Echeverr\'ia's work was partially supported by the Einstein Center for Mathematics, Berlin, the Deutscher Akademischer Austauschdienst (DAAD), Germany, and the Consejo Nacional de Ciencia y Tecnolog\'ia  (CONACyT), M{\'e}xico. We thank Michele Benzi and an anonymous referee for their    helpful comments.

\bibliographystyle{siam}
\bibliography{EchLieNab17}
\end{document}